\documentclass[10pt]{amsart}

\usepackage[margin=1.5in]{geometry}

\usepackage{amsmath, amssymb, amsthm}
\usepackage[pdftex]{color}
\usepackage{comment}
\usepackage[colorlinks,citecolor=blue,linkcolor=red]{hyperref}
\usepackage{xcolor}
\usepackage{marginnote}
\usepackage{hyperref}
\usepackage[all,arc]{xy}

\newtheorem{theorem}{Theorem}[section]
\newtheorem{lemma}[theorem]{Lemma}
\newtheorem{proposition}[theorem]{Proposition}
\newtheorem{corollary}[theorem]{Corollary}

\theoremstyle{remark}

\theoremstyle{definition}
\newtheorem{definition}[theorem]{Definition}

\newtheorem{example}{Example}

\newcommand{\define}[1]{\textbf{#1}}

\renewcommand{\P}{\mathcal{P}}
\newcommand{\PP}{\mathbb{P}}
\newcommand{\ev}{\mathrm{ev}}

\newcommand{\cay}{\mathrm{Cay}}

\newcommand{\G}{\mathfrak{G}}

\newcommand{\nc}{\newcommand}
\nc{\dmo}{\DeclareMathOperator}

\nc{\R}{\mathbb{R}}
\nc{\Z}{\mathbb{Z}}
\nc{\N}{\mathbb{N}}
\nc{\cS}{\mathcal{S}}
\nc{\iso}{\cong}
\dmo{\Diff}{Diff}
\dmo{\Homeo}{Homeo}
\dmo{\dist}{dist}
\dmo\BDiff{BDiff}
\dmo\SO{SO}
\dmo\slide{sl}
\dmo\im{im}
\dmo\id{id}
\dmo\Fix{Fix}
\dmo\Out{Out}
\dmo{\T}{\mathcal{T}}
\dmo{\Te}{\mathcal{T}^{\epsilon}}
\dmo{\Me}{\mathcal{M}^{\epsilon}}

\begin{document}

\title[Counting loxodromics]{Counting problems in graph products and relatively hyperbolic groups}

\author[I. Gekhtman]{Ilya Gekhtman}
\address{Department of Mathematics\\ 
Yale University\\ 
10 Hillhouse Ave\\ 
New Haven, CT 06520, U.S.A\\}
\email{\href{mailto:ilya.gekhtman@yale.edu}{ilya.gekhtman@yale.edu}}

\author[S.J. Taylor]{Samuel J. Taylor}
\address{Department of Mathematics\\ 
Temple University\\ 
1805 North Broad Street\ 
Philadelphia, PA 19122, U.S.A\\}
\email{\href{mailto:samuel.taylor@temple.edu}{samuel.taylor@temple.edu}}

\author[G. Tiozzo]{Giulio Tiozzo}
\address{Department of Mathematics\\ 
University of Toronto\\ 
40 St George St\\ 
Toronto, ON, Canada\\}
\email{\href{mailto:tiozzo@math.toronto.edu}{tiozzo@math.toronto.edu}}

\date{\today}
%\thanks{}

\begin{abstract}
We study properties of generic elements of groups of isometries of hyperbolic spaces.
Under general combinatorial conditions, we prove that loxodromic elements are generic (i.e. they have full density with respect to counting in balls
for the word metric) 
and translation length grows linearly.  
We provide applications to a large class of relatively hyperbolic groups and graph products, including right-angled Artin groups and right-angled Coxeter groups.
\end{abstract}

\maketitle

\section{Introduction}
Let $G$ be a finitely generated group.  One can learn a great deal about the geometric and algebraic structure of $G$ by studying its actions on various negatively curved spaces. Indeed, Gromov's theory of hyperbolic groups \cite{Gromov} provides the clearest illustration of this philosophy. However, weaker forms of negative curvature, ranging from relative hyperbolicity \cite{farb1998relatively,bowditch2012relatively,  osin2006relatively} to acylindrical hyperbolicity \cite{osin2015acylindrically,Bo}, apply to much larger classes of groups and still provide rather strong consequences.
 In all of these theories, a special role is played by the \emph{loxodromic} elements of the action, i.e. those elements which act with sink-source dynamics. In this paper, we are interested in quantifying the abundance of such isometries for the action of $G$ on a hyperbolic space $X$. We emphasize that in all but the simplest situations, the natural hyperbolic spaces that arise are not locally compact. This includes actions associated to relatively hyperbolic groups \cite{farb1998relatively}, cubulated groups \cite{KK3, hagen2014weak}, mapping class groups \cite{MM1}, and $\Out(F_n)$ \cite{BF14, handel2013free}, to name only a few. Hence, in this paper we make no assumptions of local finiteness or discreteness of the action.

Suppose that $G\curvearrowright X$ is an action by isometries on a hyperbolic space $X$. We address the question: \emph{How does a typical element of $G$ act on $X$?}

When $G$ is not amenable, the word ``typical" has no well defined meaning, and depends heavily on the averaging procedure: a family of finitely supported measures exhausting $G$. Although much is now known about measures generated from a random walk on $G$ \cite{Maher, calegari2015statistics, MaherTiozzo, mathieu2014deviation}, very little is known about counting with respect to balls in the word metric. This will be our main focus.

In more precise terms, fix a finite generating set $S$ for the group $G$. Let $B_n$ be the ball of radius $n$ about $1$ with respect to the word metric $d$ determined by $S$.
Then we call a property $P$ \emph{generic} if 
\[
\frac{\# \{g\in  B_n \ : \ g \text{ has } P \}}{ \#B_n} \to 1 \qquad \textup{as }n \to \infty .
\]
 
In this language, a refinement of our questions asks when the loxodromic elements of a particular action $G \curvearrowright X$ are generic with respect to a generating set $S$. 
It is important to note that genericity in the counting model depends on the generating set: a priori, sets may be generic with respect to one word metric, but 
not with respect to another.

The results of this paper are modeled on our previous work \cite{GTT}, where we studied the situation where $G$ is itself hyperbolic. Recall that $g \in G$ is loxodromic with respect to the action $G\curvearrowright X$ if and only if its translation length $\tau_X(g) = \lim d_X(x,g^nx) / n$ is strictly positive. In \cite{GTT}, we prove that for any isometric 
action of a hyperbolic group $G$ on a hyperbolic metric space $X$, loxodromic elements are generic, and translation length grows linearly. 
However, the genericity of loxodromic elements is in general false when the hypothesis that $G$ is hyperbolic is dropped (see Example \ref{ex:nongen}). 
In the present paper, we generalize this theorem to a much larger class of groups. Our general setup is discussed below, but here is a sample:

\begin{theorem}\label{samplethm}
Suppose that either
\begin{enumerate}
\item
$G$ is a finitely generated group which admits a geometrically finite action on a CAT($-1$) space with virtually abelian parabolic subgroups 
and $S$ an admissible generating set, or
\item 
$G$ is a right-angled Artin or Coxeter group which does not split as a direct product, and $S$ is its standard vertex generating set.
\end{enumerate} 
Then for any nonelementary isometric action $G \curvearrowright X$ on a separable hyperbolic metric space $X$
there is an $L>0$ such that 
\begin{equation}\label{genericity}
\frac{\#\{g\in B_{n}:\tau_{X}(g)\geq Ln\}}{\#B_n}\to 1.
\end{equation}
In particular, loxodromic elements are generic.
%Then the set of loxodromics of $G\curvearrowright X$ generic with respect to $S$, and in fact \eqref{genericity} holds.
\end{theorem}

In fact, our theorem applies to a more general class of relatively hyperbolic groups and graph products (see Section \ref{S:appl} for precise statements and definitions)
and in fact to any group satisfying certain combinatorial conditions. Before moving to our general framework, we state one more result which may be of independent interest. It is a direct generalization of a theorem of Gou\"e'zel, Math\'eŽus, and Maucourant \cite{gouezel2015entropy} who consider the case where $G$ is hyperbolic. 

\begin{theorem}\label{sample2}
Let $G$ be as in Theorem \ref{samplethm} with generating set $S$, and  suppose that $H$ is an infinite index subgroup of $G$.
% where either
%\begin{enumerate}
%\item
%$G$ is a finitely generated group which admits a geometrically finite action on a CAT($-1$) space with abelian parabolic subgroups 
%and $S$ an admissible generating set, or
%\item 
%$G$ is a right-angled Artin or Coxeter group which does not split as a direct product, and $S$ is its standard generating set.
%\end{enumerate} 
Then 
\[
\lim_{n \to \infty} \frac{\#(H \cap B_n)}{\#B_n} = 0.
\]
That is, the proportion of elements of $G$ of length less than $n$ which lie in $H$ goes to $0$ as $n \to \infty$.
\end{theorem}

\subsection{General framework and results} 
Our general framework is as follows.
We define a \define{graph structure} to be a pair $(G,\Gamma)$ where $G$ is a countable group and $\Gamma$ is a directed, finite graph such that: 
\begin{enumerate}
\item there is a labeled vertex $v_0$, called the \emph{initial vertex}; for every other vertex $v$ there exists a directed path from $v_0$ to $v$; 
\item every edge is labeled by a group element such that edges directed out of a fixed vertex have distinct labels. 
\end{enumerate}
By (2), there exists an \emph{evaluation map} $\textup{ev} \colon E(\Gamma) \to G$
and this map extends to the set of finite paths in $\Gamma$ by concatenating edge labels.
We denote by $\Omega_0$ the set of finite paths starting at $v_0$, 
by $S_n \subset \Omega_0$ the set of paths of length $n$, and by $\# X$ the cardinality of $X$. 
A graph structure is a \define{geodesic combing} if the evaluation map $\ev : \Omega_0 \to G$ is bijective, and each path 
in $\Omega_0$ evaluates to a geodesic in the associated Cayley graph. See Section \ref{sec:graph_struct} for details.
For each set $A \subseteq \Omega_0$ we introduce the \emph{counting measure} $P^n$ on $\Omega_0$ as 
$$P^n(A) := \frac{\#(S_n \cap A)}{\# S_n}.$$
The graph structure is \emph{almost semisimple} if the number of paths of length $n$ starting from $v_0$ has \emph{pure exponential growth}, i.e.
there exists $c >0, \lambda > 1$ such that 
$$c^{-1} \lambda^n \leq \#S_n \leq c \lambda^n$$
for each $n$. See Section \ref{sec:almost_ss} for details.

\begin{definition}
For each vertex $v$ of $\Gamma$, we denote by $L_v$ the set of loops based at $v$, and by $\Gamma_v = \ev(L_v)$ its image in $G$.
We call $\Gamma_v$ the \emph{loop semigroup} associated to $v$. 
\end{definition}

Consider an action $G \curvearrowright X$, where $X$ is a hyperbolic metric space. 
A semigroup $L < G$ is \emph{nonelementary} if it contains two independent loxodromics. 
A graph structure $(G,\Gamma)$ is \define{nonelementary} for the action $G \curvearrowright X$ if for any vertex $v$ of maximal growth (see Definition \ref{D:irred}) the loop semigroup $\Gamma_v$ is nonelementary. 

We now introduce several criteria on a graph structure that guarantee it is nonelementary: we call them thickness and quasitightness.

\begin{definition}[Thickness]
A  graph structure is \define{thick} if for any vertex $v$ of maximal growth there exists a finite set $B \subseteq G$ such that 
$$G = B \Gamma_v B.$$
More generally, a graph structure is \emph{thick relatively to a subgroup} $H < G$ if for every vertex $v$ of maximal growth 
there exists a finite set $B \subseteq G$ such that 
$$H \subseteq B \Gamma_v B.$$
\end{definition}

Given a path $\gamma$ in $\Gamma$, we say it $c$--almost contains an element $w\in G$ if $\gamma$ contains a subpath $p$ such that 
$w = a \cdot \ev(p) \cdot b$ in $G$, with $|a|, |b| \leq c$. 
We denote as $Y_{w, c}$ the set of paths in $\Gamma$ starting at the initial vertex which do not $c$--almost contain $w$. The following definition is modeled on the one found in \cite{arzhantseva2002growth}.

\begin{definition}[Growth quasitightness]
A graph structure $(G, \Gamma)$ is called \define{growth quasitight} if there exists $c > 0$ such that for every $w \in G$ the set $Y_{w, c}$ has density zero with respect to $P^n$; that is,
\[
P^n(Y_{w,c}) \to 0 \quad \text{ as } n \to 0.
\]
More generally, given a subgroup $H < G$ we say that $(G, \Gamma)$ is \emph{growth quasitight relative to }$H$ if there exists a constant $c > 0$ 
such that for every $w \in H$ the set $Y_{w, c}$ has density zero.
\end{definition}

In the most general form, the main theorem we are going to prove is the following. 

\begin{theorem}\label{T:main}
Let $G$ be a countable group of isometries of a separable, $\delta$-hyperbolic metric space $X$, and let $(G, \Gamma)$ be an almost semisimple graph structure which is either: 
\begin{enumerate}
\item nonelementary; 
\item thick relative to a nonelementary subgroup $H < G$; or 
\item growth quasitight relative to a nonelementary subgroup $H < G$.
\end{enumerate}
Then there exists $L > 0$ such that  for every $\epsilon > 0$ one has that: 
\begin{itemize}
\item[(i)]
Displacement grows linearly:
$$\frac{\#\{ g \in S_n \ : \ d_X(gx,x) \geq (L- \epsilon) n \}}{\#S_n} \to 1 \qquad \textup{as }n \to \infty.$$
\item[(ii)]
Translation length grows linearly: 
$$\frac{\#\{ g \in S_n \ : \ \tau_X(g) \geq (L- \epsilon) n \}}{\#S_n} \to 1 \qquad \textup{as }n \to \infty.$$
\item[(iii)]
As a consequence, loxodromic elements are generic: 
$$\frac{\#\{ g \in S_n \ : \ g \textup{ is }X-\textup{loxodromic} \}}{\#S_n} \to 1 \qquad \textup{as }n \to \infty.$$
\end{itemize}
\end{theorem}

If we are interested in counting with respect to balls in the Cayley graph, we get the following immediate consequence. 

\begin{corollary} \label{th:3hyp_th}
Let $G$ be a group with finite generating set $S$. Suppose that
\begin{enumerate}
\item[(i)] there is a geodesic combing for $(G, S)$;
\item[(ii)] $G$ has pure exponential growth with respect to $S$; and
\item[(iii)] the combing for $(G, S)$ satisfies at least one of the conditions (1), (2), (3) above.
\end{enumerate}
Then for any nonelementary action $G \curvearrowright X$ on a hyperbolic space, the set of loxodromic elements is generic with respect to $S$.
\end{corollary}

Note that the right-angled Artin group $G=F_{2}\times F_{3}$ of Example \ref{ex:nongen} with the standard generators has a geodesic combing and has pure exponential growth but loxodromic elements are not generic, so an additional dynamical condition (such as (1), (2), (3)) must be added. In fact, we will show that for graph products such as RAAGs and RACGs this condition amounts essentially to the group $G$ not being a product. Moreover, we will prove that the three conditions are related, namely $(3) \Rightarrow (2) \Rightarrow (1)$.

\section{Applications} \label{S:appl}

\subsection{Hyperbolic groups}

By Cannon's theorem \cite{cannon1984combinatorial}, a hyperbolic group $G$ admits a geodesic combing for \emph{any} generating set.
In fact, the language recognized by the graph is defined by choosing for each $g \in G$ the smallest word (in lexicographic order) among all words 
of minimal length which represent $g$. This is called the ShortLex representative. We proved in \cite{GTT} that this graph structure is nonelementary, hence we can apply Theorem \ref{T:main}.

\subsection{RAAGs, RACGs, and graph products}

Let $G$ be a right angled Artin or Coxeter group, and let $S$ its standard vertex generating set. A result of Hermiller and Meier \cite{HerM} implies that $(G,S)$ is ShortLex automatic. In our language, $(G,S)$ admits a geodesic combing. However, the graph $\Gamma$ parameterizing this language of geodesics does not have the correct dynamical properties needed to apply Theorem \ref{T:main}. In Section \ref{RAAGs}, we modify their construction to show that when $G$ is not a direct product, it has a graph structure with respect to the standard generators with the strongest possible dynamical properties. We then obtain the following:

\begin{theorem} \label{T:intro-RAAG}
Let $G$ be a right-angled Artin or Coxeter group which is not virtually cyclic and does not split as a product, and consider an action of $G$ on a hyperbolic, separable metric space $X$. Then eq. \eqref{genericity} holds, and loxodromic elements are generic with respect to the standard generators. 
\end{theorem}

Actually, our theorem applies to all graph products of groups with geodesic combing (Theorem \ref{th:RAAGs}). We refer the reader to Section \ref{RAAGs} for details.
Let us point out that RAAGs which are products give in fact examples of actions where loxodromics are \emph{not} generic: 
\begin{example}[Nongenericity in general] \label{ex:nongen}
Denote the free group of rank $n$ by $F_n$ and fix a free basis as a generating set. Let $G  = F_2 \times F_3$ and let $X$ denote a Cayley graph for $F_2$. Give $G$ its standard generating set; that is, the generating set consisting of a basis for $F_2$ and a basis for $F_3$.
Consider the action $G \curvearrowright X$ in which the $F_2$ factor acts by left multiplication and the right factor acts trivially. If we denote the set of loxodromics for the action by $\mathrm{LOX}$, then 
\[
\lim_{n\to \infty} \frac{\# \left(\mathrm{LOX} \cap B_n \right)}{\#B_n} = \frac{2}{3} \neq 1.
\]
Note that in the example above $G$ has pure exponential growth and a geodesic combing, so these two conditions are not sufficient to yield genericity of loxodromics. 
Moreover, the complement of $\mathrm{LOX}$ is a subgroup $H < G$ which has infinite index and positive density, showing that conditions are needed also in Theorem \ref{sample2}. 
\end{example}

\medskip
Moreover, as a consequence of the geodesic combing that we produce in order to prove the previous theorem, we also prove the following fine 
counting statement for the number of elements in a ball with respect to the standard generating set. As far as we know, this result is also new, 
and it may be of independent interest. 

\begin{theorem} \label{T:intro-exact}
Let $G$ be a right-angled Artin group or Coxeter group which is not virtually cyclic and does not split as a product. Then there exists $\lambda > 1$, $C > 0$ such that 
$$\lim_{n \to \infty} \frac{\# S_n}{\lambda^n} = C.$$
\end{theorem}

We say that a group with a generating set with the previous property has \emph{exact exponential growth}. This is stronger than \emph{pure exponential growth}
(where one only requires $C^{-1} \lambda^n \leq \# S_n \leq C \lambda^n$), 
and depends very subtly on the choice of generating set. In fact, in Theorem \ref{T:exact-exp} we will establish this result also more generally for graph products. 

\subsection{Relatively hyperbolic groups} \label{S:introRH}

Our results also apply to a large class of relatively hyperbolic groups. We need two hypotheses.
 
First, recall that a relatively hyperbolic group $G$ is equipped with a compact metric space $\partial G$ known as its \emph{Bowditch boundary}, 
and such a space carries a natural \emph{Patterson-Sullivan measure} $\nu$, defined with respect to the word metric on $\cay(G, S)$ (see Section \ref{S:PS-RH}). 
We call a relatively hyperbolic group $G$ with a generating set $S$ \emph{pleasant} if the action of $G$ on $\partial G \times \partial G$ is ergodic with respect to the measure $\nu \times \nu$. 

Second, we need a geodesic combing with respect to some generating set $S$. Let us call a finite generating set $S$ \emph{admissible} if $G$ admits a geodesic combing with respect to $S$. We have the following general statement: 

\begin{theorem} \label{T:mainRH}
Let $G$ be a relatively hyperbolic group with an admissible generating set $S$ for which $G$ is pleasant.
% such that: 
%\begin{enumerate}
%\item $G$ is pleasant; 
%\item each parabolic subgroup $P \in \mathcal{P}$ is geodesically completable.
%\end{enumerate}
Then, %each finite generating set for $G$ can be extended to a finite generating set such that 
for each action of $G$ on a hyperbolic, separable, metric space $X$, there exists $L > 0$ such that  
$$\frac{\#\{g\in B_{n}:\tau_{X}(g)\geq Ln\}}{\#B_n}\to 1 \qquad \textup{as }n \to \infty.$$
As a consequence, $X$-loxodromic elements are generic. 
\end{theorem}

In fact, by \cite{AntolinCiobanu} and \cite{Neumann-Shapiro-kleinian}, many relatively hyperbolic groups admit geodesic combings as follows.
Let us call a finitely generated group $G$ \emph{geodesically completable} if any finite generating set $S$ of $G$ can be extended to a finite generating set $S' \supseteq S$ for which there exists a geodesic biautomatic structure. Antol\'in and Ciobanu (\cite{AntolinCiobanu}, Theorem 1.5) proved that whenever $G$ is hyperbolic relative to a collection of subgroups $(P_i)$ each of which is geodesically completable, then $G$ is geodesically completable. Moreover, from automata theory (\cite{GroupsLangAuto}, Theorem 5.2.7) one gets that if $G$ admits a geodesic biautomatic structure for $S$, then it also admits a geodesic combing for the same $S$. Hence one gets: 

\begin{proposition}
Let $(G, \mathcal{P})$ be a relatively hyperbolic group such that each parabolic subgroup $P \in \mathcal{P}$ is geodesically completable. 
Then every finite generating set $S$ can be extended to a finite generating set $S'$ which admits a geodesic combing.
\end{proposition}

Let us note that in particular, virtually abelian groups are geodesically completable (\cite{AntolinCiobanu}, Proposition 10.1), hence any group hyperbolic relative to a collection of virtually abelian subgroups is geodesically completable and admits a geodesic combing. 
Moreover, we will prove (Proposition \ref{doublyergodic}): 

\begin{proposition}
If a group $G$ acts geometrically finitely on a CAT($-1$) proper metric space, then $G$ is pleasant with respect to any finite generating set. 
\end{proposition} 

In particular, geometrically finite Kleinian groups satisfy both hypotheses of Theorem \ref{T:mainRH}, 
which establishes Theorem \ref{samplethm} (1) as a corollary of Theorem \ref{T:mainRH}. 

\subsection{Actions with strongly contracting elements}

Let us now remark that by combining our work with recent work of W. Yang one can apply our theorem in even more general cases. 
Following \cite{Arzhantseva-Cashen-Tao}, \cite{Wenyuan-scc}, we call an element $g\in G$ \emph{strongly contracting} for the action on $\cay(G, S)$ in the above sense if $n \mapsto g^n$ is a quasigeodesic and there exists $C,D\ge0$ such that for any geodesics $\gamma$ in $\cay(G,S)$ whose distance from $\langle g \rangle$ is at least $C$, the diameter of the image of $\gamma$ under the nearest point projection to $\langle g \rangle$ is bounded by $D$.

Wenyuan Yang \cite{Wenyuan-scc} has recently announced that whenever the action $G\curvearrowright \cay(G,S)$ has a strongly contracting element, $G$ is growth quasitight and has pure exponential
growth with respect to $S$. 
Combining Theorem \ref{th:3hyp_th} with Yang's result we obtain the following: %corollary:

\begin{corollary} \label{C:wenyuan}
Let $G$ be a group with finite generating set $S$. Suppose that the Cayley graph $\cay(G, S)$ has a strongly contracting element and that $(G,S)$ has a geodesic combing. Then for any nonelementary action $G \curvearrowright X$ on a hyperbolic space, the set of loxodromic elements is generic with respect to $S$.
\end{corollary}

\subsection{Genericity with respect to the Markov chain}
Our approach is to deduce typical properties of elements of $G$ from typical long term behavior of long paths in the associated graph structure. 
As a by-product, we also obtain a general theorem about generic elements for sample paths in a Markov chain, which may be of independent interest. 
More precisely, an almost semisimple graph $\Gamma$ defines a Markov chain on the vertices of $\Gamma$ (see \ref{sec:Mchains}), hence it defines a Markov measure $\mathbb{P}$ on the set $\Omega_0$ of infinite paths from the initial vertex. For such Markov chains, we prove the following: 

\begin{theorem}\label{T:Pgeneric}
Let $(G, \Gamma)$ be an almost semisimple, nonelementary graph structure for $G\curvearrowright X$, and let $x \in X$. 
Then: 
\begin{enumerate}
\item 
For $\mathbb{P}$-almost every sample path $(w_n)$, the sequence $(w_n x)$ converges to a point in $\partial X$; 
\item 
There exist finitely many constants $L_i >0$ $(i = 1, \dots, r)$ such that for $\mathbb{P}$-almost every sample path 
there exists an index $i$ such that 
$$\lim_{n \to \infty} \frac{d(w_n x, x)}{n} = L_i;$$
\item
If we denote $L :=  \min_{1 \leq i \leq r} L_i$, then for each $\epsilon >0$ one has  
$$\mathbb{P}(\tau_X(w_n) \geq n (L - \epsilon)) \to 1$$
as $n \to \infty$. As a consequence, 
$$\mathbb{P}(w_n \textup{ is loxodromic}) \to 1 \quad \text{as } n\to \infty.$$
\end{enumerate} 
\end{theorem}

\subsection{Non-backtracking random walks}

An illustration of the previous result is given by looking at non-backtracking random walks. 
Let $G$ be a group and $S = S^{-1}$ a generating set. The \emph{non-backtracking random walk} on $G$ is the process defined 
by taking $g_{n}$ uniformly at random among the elements of $S \setminus \{ (g_{n-1})^{-1} \}$ and considering the sample path 
$w_n = g_1 g_2 \dots g_n$. We prove the following, which answers a question of I. Kapovich.

\begin{theorem}
Let $G$ be a nonelementary group of isometries of a separable hyperbolic metric space $X$, and let $S$ be a finite generating set. 
Consider the non-backtracking random walk  
$$w_n := g_1 \dots g_n$$
defined as above, and let $\mathbb{P}$ be corresponding the measure on the set $\Omega_0$ of sample paths. Then 
$$\mathbb{P} \left(w_n \textup{ is loxodromic on }X \right) \to 1$$
as $n \to \infty$.
\end{theorem}

\begin{proof}
Let us consider $F = F(S)$ the free group generated by $S$, with its standard word metric. By composing the projection $F \to G$ with the action on $X$, we can think of $F$ as a group of isometries of $X$. Then $F$ has a standard geodesic combing, whose graph $\Gamma$ has only one non-trivial component, hence (by Proposition \ref{P:T->NE}) the graph structure $(F, \Gamma)$ is thick, hence nonelementary. The result then follows from Theorem \ref{T:Pgeneric}.
\end{proof}

\subsection{Previous results}
Beginning with Gromov's influential works \cite{Gromov, gromov1993geometric, gromov2003random}, there 
is a large literature devoted to studying typical behavior in finitely generated groups.
More recent developments can be found, for example, in \cite{ArzhantsevaOlshanskii, arzhantseva1998generic, borovik2003multiplicative,  champetier1995proprietes, kapovich2003generic, kapovich2005genericity, kapovich2007densities, ol1992almost}. 

If one takes the definition of genericity with respect to random walks, instead of using counting in balls, then genericity of loxodromics 
%As mentioned before, genericity of loxodromics with respect to random walks, rather than counting in balls, 
has been established in many cases. 
In particular, the question of genericity of pseudo-Anosovs in the mapping class group goes back to at least Dunfield-Thurston \cite{DunfieldThurston}, 
and for random walks it has been proven independently by Rivin \cite{rivin2008walks} and Maher \cite{Maher}. 
This relates to our setup, as a mapping class is pseudo-Anosov iff it acts loxodromically on the curve complex. 
Genericity of loxodromics for random walks on groups of isometries of hyperbolic spaces has been established with increasing level of generality in 
%started with \cite{kaimanovich1994poisson} for proper actions, and led to several generalizations including
 \cite{calegari2015statistics, sisto2011contracting, MaherTiozzo}.
%In particular, Maher-Tiozzo \cite{MaherTiozzo} proved that for \emph{any} nonelementary group of isometries of a $\delta$-hyperbolic space, loxodromic elements are 
%typical along random walks. 
Let us note that in general counting in balls and counting with random walks need not yield the same result, and in fact it is a very important problem to 
establish whether the harmonic measure for the random walk can coincide with a Patterson-Sullivan-type measure, given by taking limits of counting measures over balls. 
Many results in this area show that the two measures do not coincide except in particular cases, while an existence result of a random walk for which harmonic and PS measure coincide is due for hyperbolic groups to 
%A positive answer for hyperbolic groups has been given by 
Connell-Muchnik \cite{connell2007harmonicity}.%, while many more resultsbut many negative results have been proved. 

As for counting in balls, Wiest \cite{wiest2014genericity} recently showed that if a group $G$ satisfies a weak automaticity condition and the action $G \curvearrowright X$ on a hyperbolic space $X$ satisfies a strong \emph{geodesic word hypothesis}, then the loxodromics make up a \emph{definite proportion} of elements of the $n$ ball. This geodesic word hypothesis essentially requires geodesics in the group $G$, given by the normal forms, to project to unparameterized quasigeodesics in the space $X$ under the orbit map. In our work, on the other hand, we do not assume any nice property of the action except it being by isometries. Let us note that our theorems 
answer (when the hypotheses of our two papers overlap) the open problems (3) (4) in (\cite{wiest2014genericity}, section 2.12)

\subsection*{Acknowledgments}
We thank Yago Antol\'in for useful suggestions and clarifications. The first author is partially supported by NSF grant DMS-1401875 and ERC advanced grant ``Moduli" of Prof. Ursula Hamenst\"adt, the second author is partially supported by NSF grants DMS-1400498 and DMS-1744551, and the third author is partially supported by NSERC and the Connaught fund.

\section{Background material}
Since graph structures play central role in our work, we begin by discussing some further details. The reader will notices that much of this is inspired by the theory of regular languages and automatics groups  \cite{epstein1992word}, but we place a special focus on the graph which parameterizes the language. Thus our terminology may differ from that in the literature.

\subsection{Graph structures}  \label{sec:graph_struct}
The general framework is as follows.
We define a \define{graph structure} to be a pair $(G,\Gamma)$ where $G$ is a countable group and $\Gamma$ is a directed, finite graph such that: 
\begin{enumerate}
\item there is a labeled vertex $v_0$, called the \emph{initial vertex}; for every other vertex $v$ there exists a directed path from $v_0$ to $v$; 
\item every edge is labeled by a group element such that edges directed out of a fixed vertex have distinct label. 
\end{enumerate}

Thus, there exists an \define{evaluation map} $\textup{ev} \colon E(\Gamma) \to G$
and this map extends to the set of finite paths in $\Gamma$ by concatenating edge labels. If $S = \ev (E(\Gamma))$, we say that $(G,\Gamma)$ is a graph structure with respect to $S$. 
We denote by $\Omega_0$ the set of finite paths starting at $v_0$ and by $\Omega$ the set of all finite paths.
When $\ev(\Omega_0)= G$, we call the graph structure \define{surjective}; in this case $S = \ev(E(\Gamma))$ generates $G$ as a semigroup. 
A surjective graph structure is \define{geodesic} if for each path $p \in \Omega$, the word length $|\ev(p)|_S$ is equal to the length of the path. In this case, all paths in $\Gamma$ evaluate naturally to geodesic paths in the Cayley graph $\cay(G,S)$. Finally, the graph structure is called \define{injective}, if $\ev \colon \Omega_0 \to G$ is injective. For example, if each path in $\Omega_0$ labels the ShortLex geodesic representative of its evaluation (with respect to some ordering on $S$), then $(G, \Gamma)$ is injective. A bijective, geodesic 
graph structure $(G,\Gamma)$ with respect to $S$ is called a \define{geodesic combing} of $G$ with respect to $S$.

Note the evaluation map, restricted to $\Omega_0$, factors through $S^*$, the set of all words in the alphabet $S$. The image in $S^*$ of $\Omega_0$ (i.e. all words which can be spelled starting at $v_0$) is called the language parameterized (or recognized) by $\Gamma$. This language is prefix closed by construction; an initial subword of a recognized word is also recognized. We warn the reader that references differ on the exact meaning on some of these terms. For example, Calegari--Fujiwara use the term ``combing'' to refer to the \emph{language} of a bijective, geodesic graph structure rather than the graph structure itself \cite{calegari2010combable, calegari2013ergodic}.  Since we will be most interested in dynamical properties of the graph parameterizing the language of geodesics, we choose to emphasize the graph structure itself.

\subsection{Almost semisimple graphs} \label{sec:almost_ss}
Let us summarize some of the fundamental properties about graphs and Markov chains. Much of this material appears in Calegari--Fujiwara \cite{calegari2010combable}, and we refer to that article and \cite{GTT} for more details and proofs.

Let $\Gamma$ be a finite, directed graph with vertex set $V(\Gamma) = \{v_0, v_1, \dots, v_{r-1}\}$. 
The \emph{adjacency matrix} of $\Gamma$ is the $r \times r$ matrix $M = (M_{ij})$ defined so that 
$M_{ij}$ is the number of edges from $v_i$ to $v_j$. 

Such a graph is \emph{almost semisimple} of growth $\lambda > 1$ if the following hold: 
\begin{enumerate}
\item 
There is an \emph{initial vertex}, which we denote as $v_0$; 
\item 
For any other vertex $v$, there is a (directed) path from $v_0$ to $v$; 
\item 
The largest modulus of the eigenvalues of $M$ is $\lambda$, and for any eigenvalue of modulus $\lambda$, its geometric 
multiplicity and algebraic multiplicity coincide.
\end{enumerate}

We denote by $\Omega$ the set of all finite paths in $\Gamma$, $\Omega_v$ for the set of finite paths starting at $v$, and $\Omega_0 = \Omega_{v_0}$ the set of finite paths starting at $v_0$. For a path $g \in \Omega$, we use $[g]$ to denote its terminal vertex. Similarly, we denote as $\Omega^\infty$ the set of all infinite paths in $\Gamma$, 
$\Omega_v^\infty$ the set of infinite paths starting at $v$ and $\Omega^\infty_0 = \Omega^\infty_{v_0}$.

Given two vertices $v_1$, $v_2$ of a directed graph, we say that $v_2$ is \emph{accessible from} $v_1$ and write $v_1 \to v_2$ if there is a path from $v_1$ to $v_2$, 
and two vertices are \emph{mutually accessible} if $v_1 \to v_2$ and $v_2 \to v_1$. 
Mutual accessibility is an equivalence relation, and equivalence classes are called \emph{irreducible components} of $\Gamma$. 

\medskip

For any subset $A \subseteq \Omega_0$, we define the \emph{growth} $\lambda(A)$ of $A$ as 
$$\lambda(A) := \limsup_{n \to \infty} \sqrt[n]{\# (A \cap S_n)},$$
where $S_n \subset \Omega_0$ is the set of all paths starting at $v_0$ that have length $n$.

\medskip
For each vertex $v$ of $\Gamma$ which lies in a component $C$, let $\P_v(C)$ denote the set of finite paths in $\Gamma$ based at $v$ which lie entirely in $C$. 
Moreover, for any path $g$ from $v_0$ to $v$, we let $\P_g(C) = g \cdot \P_v(C)$ be the set of finite paths in $\Omega$ which can be written as a concatenation of $g$ with a path contained entirely in $C$.

\begin{definition} \label{D:irred}
An irreducible component $C$ of $\Gamma$ is called \emph{maximal} if for some (equivalently, any) $g\in \Omega_0$ with $[g]\in C$, the growth of $\P_g(C)$ equals $\lambda$. A vertex is \emph{maximal} if it belong to a vertex of maximal growth.
Moreover, we say a vertex $v_i$ of $\Gamma$ has \emph{large growth} if there exists a path from $v_i$ to a vertex in a maximal component, and it has \emph{small growth} 
otherwise. 
\end{definition}

\begin{definition} \label{def:loopy}
For every vertex $v$ of $\Gamma$, the \emph{loop semigroup} of $v$ is the set  $L_v$ of loops in the graph 
$\Gamma$ which begin and end at $v$. It is a semigroup with respect to concatenation.
A loop in $L_v$ is \emph{primitive} if it is not the concatenation of two (non-trivial) loops in $L_v$. 
\end{definition}

\label{L:growth}
Let $\Gamma$ be an almost semisimple graph of growth $\lambda > 1$. Then there exist constants $c > 0$ and $\lambda_1 < \lambda$ such that  (\cite{GTT}, Lemma 2.3):
\begin{enumerate}
\item 
For any vertex $v$ of large growth and any $n \geq 0$, 
$$c^{-1} \lambda^n \leq \# \{ \textup{paths from }v\textup{ of length }n \} \leq c \lambda^n$$
\item 
For any vertex $v$ of small growth and any $n \geq 0$, 
$$\# \{ \textup{paths from }v\textup{ of length }n \} \leq c \lambda_1^n$$
\item 
If $v$ belongs to the maximal component $C$, then for any $n \geq 0$  
$$c^{-1} \lambda^n \leq \# \{\textup{paths in } \mathcal{P}_v(C) \textup{ of length }n \} \leq c \lambda^n $$
and also  (\cite{GTT}, Lemma 6.5)
$$c^{-1} \lambda^n \leq \# \{\textup{paths in } L_v \textup{ of length }n \} \leq c \lambda^n $$
\end{enumerate}

\subsection{Markov chains} \label{sec:Mchains}

Given an almost semisimple graph $\Gamma$ of growth $\lambda$ with edge set $E(\Gamma)$, one constructs a Markov chain on the vertices of $\Gamma$ 
as follows. If $v_i$ has large growth, then we set the probability $\mu(v_i \to v_j)$ of going from $v_i$ to $v_j$ as
\begin{equation} \label{E:Markov}
\mu(v_i \to v_j) = \frac{M_{ij} \rho(1)_j}{\lambda \rho(1)_i},
\end{equation}
and if $v_i$ has small growth, we set $\mu(v_i \to v_j) = 0$ for $i\neq j$ and $\mu(v_i \to v_i) = 1$.

\medskip 
Now, for each vertex $v$ the measure $\mu$ induces measures $\mathbb{P}_v^n$ on the space $\Omega_v$ of finite paths
starting at $v$, simply by setting 
$$\mathbb{P}_v^n(\gamma) = \mu(e_1)\cdots \mu(e_n)$$ 
for each path $\gamma = e_1 \dots e_n$ starting at $v$.
Similarly, this measure can be extended to a measure $\mathbb{P}_v$ on the space $\Omega_v^\infty$ of infinite paths starting at $v$. 
The most important cases for us will be the measures on the set of (finite and infinite, respectively) paths starting at $v_0$, which we will denote as
$\mathbb{P}^n = \mathbb{P}_{v_0}^n$ and $\mathbb{P} = \mathbb{P}_{v_0}$. 
Each measure $\mathbb{P}_v$ defines a Markov chain on the space $V(\Gamma)$, and we consider for each $n$ the random variable 
$$w_n : \Omega^\infty \to \Omega$$
$$w_n((e_1, \dots, e_n, \dots)) = e_1 \dots e_n$$
defined as the concatenation of the first $n$ edges of the infinite path.

In order to compare the $n$-step distribution for the Markov chain to the counting measure, let us denote as $LG$ the set of paths from $v_0$ ending at a vertex of large growth. Then we note (\cite{GTT}, Lemma 3.4) that 
there exists $c > 1$ such that, for each $A \subseteq \Omega_0$, 
\begin{equation} \label{E:countmu}
c^{-1} \ \mathbb{P}^n(A) \leq P^n(A \cap LG) \leq c \ \mathbb{P}^n(A).
\end{equation}

It turns out (see \cite{GTT}, Lemma 3.3) that, with respect to this choice of measure, a vertex $v$ belongs to a maximal irreducible component of $\Gamma$ if and only if it is \emph{recurrent}, i.e. :
\begin{enumerate}
\item there is a path from $v_0$ to $v$ of positive probability; and 
\item
whenever there is a path from $v$ to another vertex $w$ of positive probability, there is also a path from $w$ to $v$ 
of positive probability. 
\end{enumerate}
For this reason, maximal components will also be called \emph{recurrent components}. 

It is well-known that for almost every path of the Markov chain there exists one recurrent component $C$ such that the path 
lies completely in $C$ from some time on, and visits each vertex of $C$ infinitely many times. Thus, for each recurrent component $C$, we let $\Omega_C$ be the set of all infinite paths from the initial vertex which enter $C$ and remain inside $C$ forever, and denote as $\mathbb{P}_C$ the conditional probability of $\mathbb{P}$ on $\Omega_C$. 

Moreover, for each recurrent vertex $v$ 
the distribution of return times decays exponentially: 
\begin{equation} \label{E:return}
\mathbb{P}_v\left( \tau^+_v = n \right) \leq e^{-cn}
\end{equation}
where $\tau^+_v = \min \{ n \geq 1 \ : \ [w_n] = v \}$ denotes the first return time to vertex $v$. 

We will associate to each recurrent vertex of the Markov chain a random walk, and use previous results on random walks to 
prove statements about the asymptotic behavior of the Markov chain. 

\medskip
For each sample path $\omega \in \Omega^\infty$, let us define $n(k,v, \omega)$ as the $k^{th}$ time the path $\omega$ lies at the vertex $v$. 
In formulas, 
$$n(k, v, \omega) := \left\{ \begin{array}{ll} 0 & \textup{if } k = 0 \\ \min \left( h > n(k-1, v, \omega) \ : \ [w_h] = v \right) &  \textup{if }k \geq 1 \end{array} \right.$$
To simplify notation, we will write $n(k,v)$ instead of $n(k,v, \omega)$ when the sample path $\omega$ is fixed.  

\medskip

We now define the \emph{first return measure} $\mu_v$ on the set %$L_v$ 
of primitive loops by setting, for each primitive loop $\gamma = e_1 \dots e_n$ 
with edges $e_1, \dots, e_n$, 
$$\mu_v(e_1 \dots e_n) = \mu(e_1)\dots \mu(e_n).$$
Extend $\mu_v$ to the entire loop semigroup $L_v$ by setting $\mu_v(\gamma) = 0$ if $\gamma\in L_v$ is not primitive. 
Since almost every path starting at $v$ visits $v$ infinitely many times, the measure $\mu_v$ is a probability measure.

By equation \eqref{E:return}, for every recurrent vertex $v$, the first return measure $\mu_v$ has finite exponential moment, 
i.e. there exists a constant $\alpha >0$ such that 
\begin{equation} \label{e:expo}
\int_{L_v} e^{\alpha|\gamma|} \ d\mu_v(\gamma) < \infty.
\end{equation}

\subsection{Hyperbolic spaces} \label{sec:hyperbolic}

In this paper, $X$ will always be a geodesic, separable metric space. Such a space is called \emph{$\delta$-hyperbolic} for some $\delta \ge 0$
if for every geodesic triangle in $X$, each side is contained within the $\delta$--neighborhood of the other two sides. Each hyperbolic space has a well-defined \emph{Gromov boundary} $\partial X$, and we refer the reader to \cite[Section III.H.3]{BH}, \cite{GhysdelaHarpe}, or \cite[Section 2]{kap_boundaries} for definitions and properties. 

If $g$ is an isometry of $X$, its \emph{translation length} is defined as 
$$\tau_X(g) := \lim_{n \to \infty} \frac{d(g^n x, x)}{n}$$
where the limit does not depend on the choice of $x$. 
In order to estimate the translation length, we will use the following well-known lemma; see for example \cite[Proposition 5.8]{MaherTiozzo}.

\begin{lemma}\label{l:tau-formula}
There exists a constant $c$, which depends only on $\delta$, such that for any isometry $g$ of a $\delta$-hyperbolic space $X$ 
and any $x \in X$ with $d(x, gx) \geq 2(gx, g^{-1}x)_x + c$, the translation length of $g$ is given by 
$$\tau_X(g) = d(x, gx) - 2(gx, g^{-1}x)_x + O(\delta).$$ %\colorbox{pink}{still true}
\end{lemma}

An isometry $g$ of $X$ is \emph{loxodromic} if it has positive translation length; in that case, it has two fixed points on $\partial X$. 
We say two loxodromic elements are \emph{independent} if their fixed point sets are disjoint. A semigroup (or a group) $G < \textup{Isom }X$ is \emph{nonelementary} if it contains two independent loxodromics. 
We will use the following criterion. 

\begin{proposition}[{\cite[Proposition 7.3.1]{das2014geometry}}] \label{th:Tushar}
Let $L$ be a semigroup of isometries of a hyperbolic metric space $X$. If the limit set $\Lambda_L \subset \partial X$ of $L$ on the boundary of X is nonempty and $L$ does not have a finite orbit in $\partial X$, then $L$ is nonelementary. 
\end{proposition}

Finally, we turn to the definition and basic properties of shadows in the $\delta$-hyperbolic space $X$. For $x,y \in X$, the \emph{shadow in $X$ around $y$ based at $x$} is 
\[
S_x(y,R) = \{z\in X : (y,z)_x \ge d(x,z) -R\},
\]
where $R> 0$ and $(y,z)_x = \frac{1}{2}(d(x,y) +d(x,z) -d(y,z))$ is the usual \emph{Gromov product}. The \emph{distance parameter} of $S_x(y,R)$ is by definition the number $r = d(x,y) -R$, which up to an additive constant depending only on $\delta$, measures the distance from $x$ to $S_x(y,R)$. Indeed, $z \in S_x(y,R)$ if and only if any geodesic $[x,z]$  $2\delta$--fellow travels any geodesic $[x,y]$ for distance $r+O(\delta)$. 
The following observation is well-known.

\begin{lemma} \label{lem:neigh-shadow}
For each $D \geq 0$, and each $x, y$ in a metric space, we have 
$$N_D(S_x(y, R)) \subseteq S_x(y, R + D).$$
\end{lemma}

\subsection{Random walks}
A probability measure $\mu$ on $G$ is said to be \emph{nonelementary} with respect to the action $G \curvearrowright X$ 
if the semigroup generated by the support of $\mu$ is nonelementary.

We will need the fact that a random walk on $G$ whose increments are distributed according to a nonelementary measure $\mu$ almost surely converge to the boundary of $X$ and has positive drift in $X$.

\begin{theorem}[{\cite[Theorems 1.1, 1.2]{MaherTiozzo}}] \label{thm:MT_drift}
Let $G$ be a countable group which acts by isometries on a separable hyperbolic space $X$, and let $\mu$ be a nonelementary probability distribution on $G$. Fix $x\in X$, and let $(u_n)$ be the sample path of a random walk with independent increments with distribution $\mu$. Then:
\begin{enumerate}
\item
almost every sample path $(u_nx)$ converges to a point in the boundary of $\partial X$, and the resulting hitting measure $\nu$ is nonatomic; 
\item
moreover, if $\mu$ has finite first moment, then there is a constant $L>0$ such that for almost every sample path
$$\lim_{n \to \infty}\frac{d(x,u_nx)}{n} = L >0.$$
\end{enumerate}
\end{theorem}
\noindent The constant $L >0$ in Theorem \ref{thm:MT_drift} is called the \emph{drift} of the random walk $(u_n)$. 

\section{Behavior of generic sample paths for the Markov chain}
\label{sec:Markov_action}

Let $G$ be a group with a nonelementary action $G \curvearrowright X$ on a hyperbolic space $X$. 
In this section we assume that $G$ has a graph structure $(G,\Gamma)$ which is almost semisimple and nonelementary. 

\subsection{Convergence to the boundary of $X$} \label{sec:convergence}

In this section we show that almost every sample path for the Markov chain converges to the boundary of $X$. 
Since we are \emph{assuming} that the graph structure is nonelementary, the exact same proof as in (\cite{GTT}, Theorem 6.8) yields the following.

\begin{theorem} \label{th:chain_converges}
For $\mathbb{P}$-almost every path $(w_n)$ in the Markov chain, the projection $(w_n x)$ to the space converges to a point in the boundary $\partial X$. 
\end{theorem}

As a consequence, we have for every vertex $v$ of large growth a well-defined \emph{harmonic measure} $\nu_v^X$, namely the hitting measure 
for the Markov chain on $\partial X$: for each each (Borel) $A \subseteq \partial X$ we define 
$$\nu_v^X(A) := \mathbb{P}_v(\lim_{n \to \infty} w_n x \in A).$$

The previous proof also provides a decomposition theorem for the harmonic measures $\nu_v^X$. 
If $\mathcal{R}$ is the set of recurrent vertices of $\Gamma$, then we have:
\begin{equation} \label{E:combo}
\nu_v^X = \sum_{w \in \mathcal{R}} \sum_{\gamma: v \to w} \mu(\gamma) \gamma_* \nu_w^X
\end{equation}
Here, the sum is over all finite paths from $v$ to $w$ which only meet a recurrent vertex at their terminal endpoint.
Note that if $v$ is recurrent, then $\nu_v^X$ is the harmonic measure for the random walk on $G$ generated by the measure $\mu_v$, as discussed above. 

\begin{lemma} \label{lem:nonatomic}
For any $v$ of large growth, the measure $\nu_v^X$ is non-atomic.
\end{lemma}

\begin{proof}
Since the random walk measures $\nu^X_w$ are non-atomic, so are the measures $\gamma_* \nu^X_w$ for each $\gamma$, 
hence by equation \eqref{E:combo} the measure $\nu_v^X$ is also non-atomic as it is a linear combination of non-atomic measures. 
\end{proof}

\subsection{Positive drift along geodesics} \label{sec:drift}

In this section, we show that almost every sample path has a well-defined and positive drift in $X$. 

\begin{theorem} \label{th:drift}
For $\mathbb{P}$-almost every sample path $\omega = (w_n)$ there exists a recurrent component $C = C(\omega)$ for which we have 
$$\lim_{n \to \infty} \frac{d(w_n x, x)}{n} = L_C,$$
where $L_C >0$ depends only on $C$.
\end{theorem}

Since $\Gamma$ is finite, this gives at most finitely many potential drifts for the Markov chain.

\begin{proof}
Let $v$ be a recurrent vertex. Since the graph structure is nonelementary, the loop semigroup $\Gamma_v$ is nonelementary, hence the random walk 
given by the return times to $v$ has positive drift. More precisely, 
from Theorem \ref{thm:MT_drift}, there exists a constant $\ell_v > 0$ such that 
for almost every sample path which enters $v$,
$$\lim_{k \to \infty} \frac{d(w_{n(k,v)} x, x)}{k} = \ell_v.$$
Morever, as the distribution of return times has finite exponential moment, for almost every sample path one has 
$$\lim_{k \to \infty} \frac{n(k,v,\omega)}{k} = T_v.$$
These two facts imply
$$\lim_{k \to \infty} \frac{d(w_{n(k,v)} x, x)}{n(k,v)} = \frac{\ell_v}{T_v}.$$
Now, almost every infinite path visits every vertex of some recurrent component infinitely often.
Thus, for each recurrent vertex $v_i$ which belongs to a component $C$, there exists a constant $L_i > 0$ such that for $\mathbb{P}_C$-almost every path $(w_n)$, there is a limit 
$$L_i = \lim_{k \to \infty} \frac{d(w_{n(k, v_i)} x, x)}{n(k, v_i)}.$$
Let $C$ be a maximal component, and $v_1, \dots, v_k$ its vertices. Our goal now is to prove that $L_1 = L_2 = \dots L_k$. 
Let us pick a path $\omega \in \Omega_0$ such that the limit $L_i$ above exists for each $i = 1, \dots, k$, and define $A_i = \{ n(k, v_i), k \in \mathbb{N} \}$, and the equivalence relation $i \sim j$ if $L_{i} = L_{j}$. 
Since $w_{n(k, v_i)}$ and $w_{n(k, v_i) + 1}$ differ by one generator, $d(w_{n(k, v_i)}x, w_{n(k, v_i) + 1}x)$ is uniformly bounded, hence
$$\lim_{k \to \infty} \frac{d(w_{n(k, v_i) + 1} x, x)}{n(k, v_i) + 1} = \lim_{k \to \infty} \frac{d(w_{n(k, v_i)}x, x)}{n(k, v_i)} = L_i$$
so the equivalence relation satisfies the hypothesis of (\cite{GTT}, Lemma 6.9), 
hence there is a unique limit $L_C = L_i$ so that 
$$\lim_{n \to \infty} \frac{d(w_n x, x)}{n} = L_C.$$
\end{proof}

\begin{corollary} \label{cor:drift2}
For every vertex $v$ of large growth, and for $\mathbb{P}_v$-almost every sample path $(w_n)$ there exists a recurrent component $C$ accessible from $v$ such that 
we have 
$$\lim_{n \to \infty} \frac{d(w_n x, x)}{n} = L_C.$$ 
\end{corollary}

\begin{proof}
By Theorem \ref{th:drift}, for $\mathbb{P}$-almost every path which passes through $v$, 
the drift equals $L_C$ for some recurrent component $C$. 
Let $g_0$ be a path from $v_0$ to $v$ of positive probability. Then for any $\omega = (w_n) \in \Omega_v^\infty$, the path 
$(\widetilde{w}_{n}) = g_0 \cdot \omega$ belongs to $\Omega_0^\infty$, and moreover  $\widetilde{w}_{n+k} = w_{n+k}(g_0 \cdot \omega) = g_0 \cdot w_{n}(\omega)$ where 
$k = |g_0|$. Hence
\[
d(w_{n}x,x) - d(x,g_0x) \le d(\widetilde{w}_{n+k}x , x) = d(w_{n}x, g_0 x) \le d(w_{n}x,x) + d(x,g_0x)
\]
and so by Theorem \ref{th:drift}
$$\lim_{n \to \infty} \frac{d(w_n x, x)}{n} = \lim_{n \to \infty} \frac{d(\widetilde{w}_{n+k} x, x)}{n} = L_C$$
as required.
\end{proof}

For application in Section \ref{sec:counting}, we will need the following convergence in measure statement. Let us denote 
\begin{equation} \label{E:small-drift}
L := \min_{C \text{ recurrent}} L_C > 0
\end{equation}
the smallest drift.

\begin{corollary} \label{cor:drift_in_measure}
For any $\epsilon > 0$, and for any $v$ of large growth, 
$$\mathbb{P}_v\left(\frac{d(x, w_n x)}{n} \leq L - \epsilon \right) \to 0.$$
\end{corollary}

\begin{proof}
By the theorem, the sequence of random variables $X_n = \frac{d(w_nx, x)}{n}$ converges almost surely to a function $X_\infty$ with the finitely many values 
$L_1, \dots, L_r$. Moreover, for every $n$ the variable $X_n$ is bounded above by the Lipschitz constant of the orbit map $G \to X$.
Thus, $X_n$ converges to $X_\infty$ in $L^1$, yielding the claim. 
\end{proof}

\subsection{Decay of shadows for $\mathbb{P}$}

For any shadow $S$, we denote its closure in $X \cup \partial X$ by $\overline S$. Since the harmonic measures $\nu_v^X$ for the Markov chain are nonatomic (by Lemma \ref{lem:nonatomic}), we get by the same proof as in \cite{GTT} the following decay of shadows results. 

\begin{proposition}[Decay of shadows, \cite{GTT}, Proposition 6.18] \label{P:decay1}
There exists a function $f : \mathbb{R} \to \mathbb{R}$ such that $f(r) \to 0$ as $r \to \infty$, and such that 
for each $v$ of large growth, 
$$\nu_v^X \big(\overline{S_x(gx, R)}\big) \leq f(r),$$
where $r = d(x, gx) -R$ is the distance parameter of the shadow. 
\end{proposition}

\begin{proposition}[\cite{GTT}, Proposition 6.19]  \label{p:decay}
There exists a function $p : \mathbb{R}^+ \to \mathbb{R}^+$ with $p(r) \to 0$ as $r \to \infty$, such that for each vertex $v$
and any shadow $S_x(gx, R)$ we have 
\[
\mathbb{P}_v \Big( \exists n  \geq  0 \ : \ w_n x \in S_x(gx, R) \Big)  \leq p(r),
\]
where $r = d(x, gx) - R$ is the distance parameter of the shadow. 
\end{proposition}

\section{Generic elements with respect to the counting measure} \label{sec:counting}

We now use the results about generic paths in the Markov chain to obtain results about generic paths with respect to the counting measure. 

\subsection{Genericity of positive drift}
The first result is that the drift is positive along generic paths: 

\begin{theorem} \label{th:gen_drift}
Let $(G, \Gamma)$ be an almost semisimple, nonelementary graph structure, and $L > 0$ be the smallest drift as given by eq. \eqref{E:small-drift}. 
Then for every $\epsilon > 0$ one has 
$$\frac{\#\{ g \in S_n \ : \ d(gx,x) \geq (L- \epsilon) \ |g| \}}{\#S_n} \to 1 \qquad \textup{as }n \to \infty.$$
\end{theorem} 
The result follows from Corollary \ref{cor:drift_in_measure} similarly as in the proof of (\cite{GTT}, Theorem 5.1). 

\begin{proof}
Let $A_L$ denote the set of paths
$$A_L :=  \{ g \in \Omega_0 \ : \ d(gx,x) \leq L |g| \}. $$
We know by Corollary \ref{cor:drift_in_measure} that for any $L' < L$ one has 
$$\mathbb{P}^n(A_{L'}) \to 0 \qquad \textup{as }n \to \infty.$$
Then, for each path $g$ of length $n$ we denote as $\widehat{g}$ the prefix of $g$ of length $n - \lfloor \log n \rfloor$, and we observe that 
$$P^n(A_{L-\epsilon}) \leq \frac{\# \{ g \in S_n \ : \ \widehat{g} \notin LG \}}{\#S_n} + \frac{\#\{ g \in S_n \cap A_{L-\epsilon} \ :  \ \widehat{g} \in LG \}}{\#S_n}$$
and we know by \cite[Proposition 2.5]{GTT} that the first term tends to $0$. 
Now, by writing
 $g = \widehat{g}h$ with $| h | = \log |g|$ we have that $d(gx, x) \leq (L-\epsilon) |g|$ implies 
$$d(\widehat{g}x, x) \leq d(gx,x) + d(\widehat{g}x, gx) \leq (L-\epsilon) |g| + d(x, hx) \leq $$
hence, there exists $C$ such that it is less than 
$$ \leq (L-\epsilon) |g| + C \log |g| \leq L'  |\widehat{g}|$$
for any $L-\epsilon < L' < L$ whenever $|g|$ is sufficiently large. This proves the inclusion
$$\{ g \in S_n \cap A_{L-\epsilon} \ :  \ \widehat{g} \in LG \} \subseteq \{ g \in S_n \ :  \ \widehat{g} \in A_{L'} \cap LG \}$$
and by Lemma \ref{L:growth} (1) 
$$\#\{ g \in S_n \ :  \ \widehat{g} \in A_{L'} \cap LG \} \leq c \lambda^{\log n} \#( S_{n - \log n} \cap A_{L'} \cap LG) \leq $$
hence by equation \eqref{E:countmu} and considering the size of $S_{n-\log n}$
$$\leq c_1 \lambda^{\log n} \mathbb{P}^{n-\log n}(A_{L'}) \# S_{n-\log n} \leq c_2 \lambda^n \mathbb{P}^{n-\log n}(A_{L'}).$$
Finally, using that $\mathbb{P}^{n-\log n}(A_{L'}) \to 0$ we get
$$\limsup_{n \to \infty} \frac{\#\{ g \in S_n \cap A_{L-\epsilon} \ :  \ \widehat{g} \in LG \}}{\#S_n} \leq \limsup_{n \to \infty} c_3 \mathbb{P}^{n-\log n}(A_{L'})  = 0$$
which proves the claim.
\end{proof}

\subsection{Decay of shadows for the counting measure}
For $g \in G$, we set
 \[
 S_x^\Gamma(gx,R) = \left \{ h\in \Omega_0 : hx \in  S_x (gx,R) \right \}, 
 \]
where as usual, $S_x (gx,R)$ is the shadow in $X$ around $gx$ centered at the basepoint $x\in X$ and $hx = \ev(h)x$.  We will need the following decay property for $S_x^\Gamma(gx,R) \subseteq \Omega_0$.

\begin{proposition} \label{P:counting-decay}
There is a function $\rho : \mathbb{R}^+ \to \mathbb{R}^+$ with $\rho(r) \to 0$ as $r \to \infty$ such that for every $n \geq 0$
\[
P^n(S^\Gamma_x(gx, R)) \le \rho \left(d(x,gx) -R \right). 
\]
\end{proposition}

We start with the following lemma in basic calculus. 
\begin{lemma} \label{L:ptilde}
Let $p : \mathbb{R}^+ \to [0,1]$ be a decreasing function with $p(x) \to 0$ as $x \to \infty$.
For each $\alpha > 0$ and each $C > 0$, there exists a function $\widetilde{p} : \mathbb{R}^+ \to \mathbb{R}^+$ such that: 
\begin{enumerate}
\item $\widetilde{p}(x) \geq p(x)$ for each $x \in \mathbb{R}^+$; 
\item $\widetilde{p}(x+mC) \geq \widetilde{p}(x) e^{-\alpha m}$ for each $x \in \mathbb{R}^+$, $m \in \mathbb{N}$: 
\item $\widetilde{p}(x) \to 0$ as $x \to \infty$.
\end{enumerate}
\end{lemma}

\begin{proof}
Let us define $\widetilde{p}$ inductively as follows: 
$$\widetilde{p}(x) := p(x) \qquad \textup{if }x \in [0, C)$$
$$\widetilde{p}(x + C) := \max \{ p(x+C), \widetilde{p}(x) e^{-\alpha} \}$$
Thus (1) and (2) are immediate from the definition. To prove (3), note that for each $n$ there exists $k \leq n$ such that 
$\widetilde{p}(x + nC) = p(x + (n-k)C)e^{-\alpha k}$. Hence as $x + nC \to \infty$ either $x+(n-k)C \to \infty$ or $k \to \infty$, thus the only limit point is $0$.
\end{proof}

\begin{proof}[Proof of Proposition \ref{P:counting-decay}] 
Pick a path $h \in \Omega_0$ of length $n$ in $S_x^\Gamma(gx, R)$, and let $\widehat{h}$ denote the longest subpath of $h$ starting at the initial vertex and which ends in a vertex of large growth. 
Let us write $h = \widehat{h}l$ where $l$ is the second part of the path. Note that we have
$$d(hx, \widehat{h}x) = d(\widehat{h}l x, \widehat{h}x) = d(lx, x) \leq k C$$
where $k := |l|$ and $C$ is the Lipschitz constant of the orbit map, hence by Lemma \ref{lem:neigh-shadow}
$$\widehat{h} \in S_x^\Gamma(gx, R' + kC)$$
where $R' = R + D$ and $D = O(\delta)$. Note that for each element $\widehat{h}$ there are at most $c \lambda_1^k$ choices of the continuation $l$, hence 
$$\#( S_n \cap S^\Gamma_x(gx, R) ) \leq c \sum_{k = 0}^n \lambda_1^k \#(S_{n-k} \cap S_x^\Gamma(gx, R' + kC) \cap LG) \leq$$
and by using eq. \eqref{E:countmu} and Proposition \ref{p:decay}
$$\leq c_1 \sum_{k =0}^n \lambda_1^k \lambda^{n-k}\ \mathbb{P}^{n-k}(S_x^\Gamma(gx, R' + kC)) \leq 
c_1 \lambda^n \sum_{k =0}^n (\lambda_1/\lambda)^k\  p(d(x, gx) - R' - kC) $$
Now, by Lemma \ref{L:ptilde} we can replace $p$ by $\widetilde{p}$, choosing $\alpha$ so that $e^{\alpha} \lambda_1 < \lambda$, thus getting 
$$\widetilde{p}(d(x, gx) - R' - kC)  \leq e^{\alpha k} \widetilde{p}(d(x, gx) - R')$$ 
Thus, the previous estimate becomes
$$P^n( S^\Gamma_x(gx, R) ) \leq c_2 \sum_{k =0}^n (e^\alpha \lambda_1 /\lambda)^k \ \widetilde{p}(d(x, gx) - R') \leq c_3 \  \widetilde{p}(d(x, gx) - R')$$
which proves the lemma if one sets $\rho(r) := c_3 \ \widetilde{p}(r - D)$.
\end{proof}

\subsection{Genericity of loxodromics}

We now use the previous counting results to prove that loxodromic elements are generic with respect to the counting measure. 

The strategy is to apply the formula of Lemma \ref{l:tau-formula} to show that translation length grows linearly as function of the length of the path: in order to do so, one needs to show that the distance $d(gx, x)$ is large (as we did in Theorem \ref{th:gen_drift}) 
and, on the other hand, the Gromov product $(gx, g^{-1}x)_x$ is not too large. The trick to do this is to split the path $g$ in two subpaths of roughly the same length, 
and show that the first and second half of the paths are almost independent.

To define this precisely, for each $n$ let us denote $n_1 = \lfloor \frac{n}{2} \rfloor$ and $n_2 = n - n_1$. 
For each path $g \in \Omega$, we define its \emph{initial part} $i(g)$ to be the subpath given by the first $n_1$ edges of $g$, and its \emph{terminal part}
$t(g)$ to be the subpath given by the last $n_2$ edges of $g$. With this definition, $g = i(g) \cdot t(g)$ and $|i(g)| = n_1$, $|t(g)| = n_2$.
Moreover, we define the random variables $i_n, t_n : \Omega^\infty \to \Omega$ by $i_n(w) = i(w_n)$ and $t_n(w) = t(w_n)$. 
Note that by definition $i_n = w_{n_1}$ and by the Markov property we have for each paths $g, h \in \Omega$:
$$\mathbb{P}(i_n = g \textup{ and }t_n =h) = \mathbb{P}(w_{n_1} = g)\mathbb{P}_v(w_{n_2} = h)$$
where $v = [g]$. In the next Lemma, we use the notation $C(\omega)$ to refer to the recurrent component to which the sample path $\omega = (w_n)$
eventually belongs, as in Theorem \ref{th:drift}.

\begin{lemma} \label{L:drift-un}
For any $\epsilon > 0$ we have 
$$\mathbb{P}\left( \left|\lim_{n \to \infty} \frac{d(x, t_n(\omega) x)}{n} - \frac{L_{C(\omega)}}{2}\right| > \epsilon \right) \to 0$$
as $n \to \infty$.
\end{lemma}

\begin{proof}
Note that by definition $t_n(\omega) = w_{n_2}(T^{n_1} \omega)$, where $T : \Omega^\infty \to \Omega^\infty$ is the shift in the space of infinite paths.
Note that for every $A \subseteq \Omega^\infty$ by the Markov property we have 
$$\mathbb{P}(T^{-n}A) = \sum_{v \in V} \mathbb{P}([w_n] = v) \mathbb{P}_v(A)$$
Let us define the function 
$$S_n(\omega, \omega') := \left| \frac{d(x, w_n(\omega) x)}{n} - L_C(\omega') \right|. $$
Note that from Corollary \ref{cor:drift2} for every vertex $v$ of large growth and every $\epsilon > 0$
$$\mathbb{P}_v(S_n(\omega, \omega) \geq \epsilon) \to 0$$ 
Moreover, for every $n$, if the path $(e_1, \dots, e_n, \dots)$ lies entirely in the component $C$ from some point on, then the same is true for 
the shifted path $(e_{n+1}, e_{n + 2}, \dots)$, i.e. $C(T^n \omega) = C(\omega)$ almost surely, and so 
$$S_n(\omega, T^k \omega') = S_n(\omega, \omega') \qquad \textup{for all }n, k$$
hence 
$$\mathbb{P}(S_{n_2}(T^{n_1} \omega, \omega) > \epsilon) = \mathbb{P}(S_{n_2}(T^{n_1} \omega, T^{n_1}\omega) > \epsilon) = $$
$$ = \sum_{v \in V} \mathbb{P}([w_{n_1}] = v]) \mathbb{P}(S_{n_2}(\omega, \omega) > \epsilon) \leq \sum_{v \in V \cap LG} \mathbb{P}_v(S_{n_2}(\omega, \omega) > \epsilon)$$
and the right-hand side tends to $0$ by Corollary \ref{cor:drift2}, proving the claim. 
\end{proof}

We now show that $i(g)$ and $t(g)^{-1}$ generically do not fellow travel. For the argument, let $S_n(v)$ denote the set of paths in $\Omega$ which start at $v$ and have length $n$.

\begin{lemma} \label{l:prod0}
Let $f : \mathbb{R} \to \mathbb{R}$ be any function such that $f(n) \to +\infty$ as $n \to +\infty$. Then 
\[
P^n \Big( g\in \Omega_0 \ : \ (i(g) x, t(g)^{-1} x)_x \geq f(n) \Big) \to 0
\]
as $n \to \infty$. 
\end{lemma}

\begin{proof}
We compute
$$P^n \left(g \ : \ (i(g) x, t(g)^{-1}x)_x \geq f(n)\right) = \frac{\#\{ g \in S_{n_1}, [g] = v, h \in S_{n_2}(v)  \ : \ \ (gx, h^{-1}x)_x \geq f(n) \}}{\#S_n} \leq $$
and by fixing $v$ and forgetting the requirement that $[g] = v$ we have
$$\leq  \sum_{v \in V} \frac{\#\{ g \in S_{n_1}, h \in S_{n_2}(v) \ :  \ (gx, h^{-1}x)_x \geq f(n) \}}{\#S_n} \leq $$
then by fixing a value of $h$
$$\leq \sum_{v \in V} \frac{1}{\#S_n} \sum_{h \in S_{n_2}(v)} \#\left\{ g \in S_{n_1}  \ :  \ gx \in S_x(h^{-1}x, d(x, h^{-1}x) - f(n)) \right \} \leq$$
hence from decay of shadows (Proposition \ref{P:counting-decay}) follows that 
$$\leq \sum_{v \in V} \sum_{h \in S_{n_2}(v)} \frac{ \rho(f(n)) \#S_{n_1} }{ \#S_n } \leq  \frac{ \#V \#S_{n_1} \# S_{n_2} \rho(f(n))}{\#S_n} \leq c \rho(f(n)) \to 0 .$$
\end{proof}

Once we have shown that $i(g)$ and $t(g)^{-1}$ are almost independent, we still need to show that also $g$ and $g^{-1}$ are almost independent. 
In order to do so, we note that $i(g)$ is the beginning of $g$ while $t(g)^{-1}$ is the beginning of $g^{-1}$, and then we use the following 
trick from hyperbolic geometry. See e.g. \cite{TT}.

\begin{lemma}[Fellow traveling is contagious]\label{l:fellow_travel}
Let $X$ be a $\delta$--hyperbolic space with basepoint $x$ and let that $A\ge0$. If $a,b,c,d$ are points of $X$ with $(a\cdot b)_{x} \ge A$, $(c\cdot d)_{x} \ge A$, and $(a\cdot c)_{x} \le A-3\delta$. Then $(b \cdot d)_{x} -2\delta \le (a \cdot c)_{x} \le (b \cdot d)_{x} +2\delta$.
\end{lemma}

In order to apply Lemma \ref{l:fellow_travel}, we need to check that the first half of $g$ (which is $i(g)$) and the first half of $g^{-1}$ (which is $t(g)^{-1}$) 
generically do not fellow travel. 

\begin{lemma} \label{l:prod1}
For each $\eta > 0$, the probability
\[
P^n\left(g \in \Omega_0 \ : \ (t(g)^{-1} x, g^{-1} x)_x \leq \frac{n (L - \eta)}{2} \right) \to 0
\]
and 
\[
P^n\left(g \in \Omega_0 \ : \ (i(g) x, g x)_x \leq \frac{n (L - \eta)}{2}  \right) \to 0
\]
as $n \to \infty$. 
\end{lemma}

\begin{proof}
Consider the set $B_{L} := \{ g \in \Omega_0 \ : \ d(x, gx) - d(x, i(g)x) \leq \frac{|g| L}{2} \}$. 
We know by Theorem \ref{th:drift} for $\mathbb{P}$-almost every sample path we have 
$$\lim_{n \to \infty} \frac{d(x, w_n x) - d(x, w_{n_1} x)}{n} = \frac{L_C}{2} \geq \frac{L}{2}$$
hence for any $L' < L$ one has $\mathbb{P}^n(B_{L'}) \to 0$ as $n \to \infty$.
Hence, as in the proof of Theorem \ref{th:gen_drift} we get for any $\epsilon > 0$, 
\begin{equation} \label{E:g-a}
P^n\left(d(x,gx) - d(x, i(g)x) \geq \frac{n(L-\epsilon)}{2}\right)\to 1
\end{equation}
Finally, by writing out the Gromov product, the triangle inequality and the fact that the action is isometric we get 
$$(t(g)^{-1}x, g^{-1}x)_x \geq d(x, g^{-1}x) - d(t(g)^{-1} x, g^{-1} x) = d(x, gx) - d(x, i(g)x)$$
which combined with \eqref{E:g-a} proves the first half of the claim. 

The second claim follows analogously. Namely, from Theorem \ref{th:drift} and Lemma \ref{L:drift-un}, we have for any $\epsilon > 0$
$$\mathbb{P}\left( d(x, w_n(\omega) x) - d(x, t_n(\omega) x) \leq \frac{n(L_{C(\omega)} - \epsilon)}{2} \right) \to 0$$
which then implies as before  
$$P^n\left(d(x,gx) - d(x, t(g)x) \geq \frac{n(L-\epsilon)}{2}\right)\to 1$$
and to conclude we use that 
$$(i(g)x, gx)_x \geq d(x, gx) - d(i(g)x, gx) = d(x, gx) - d(x, t(g)x).$$
\end{proof}

We now use Lemma \ref{l:fellow_travel} (fellow traveling is contagious) to show that the Gromov products $(gx,g^{-1}x)_x$ do not grow too fast with respect to our counting measures. 

\begin{proposition} \label{p:g-prod}
Let $f : \mathbb{N} \to \mathbb{R}$ be a function such that $f(n) \to +\infty$ as $n \to \infty$.
Then 
$$P^n \Big( (gx, g^{-1}x)_x \leq f(n) \Big) \to 1$$
as $n \to \infty$. 
\end{proposition}

\begin{proof}
Define 
$$f_1(n) = \min \left\{ f(n) - 2 \delta, \frac{n(L-\eta)}{2} - 3\delta \right\}$$
It is easy to see that $f_1(n) \to \infty$ as $n \to \infty$.
By Lemma \ref{l:fellow_travel}, if we know that:
\begin{enumerate}
\item
 $(i(g)x, gx)_x \geq n(L-\eta)/2$, 
 \item
 $(t(g)^{-1}x, g^{-1}x)_x \geq n(L-\eta)/2$, and
\item
 $(i(g)x, t(g)^{-1}x)_x \leq f_1(n) \leq n(L-\eta)/2 - 3\delta$, 
\end{enumerate}
then 
$$(gx, g^{-1}x)_x \leq (i(g)x, t(g)^{-1}x)_x + 2 \delta \leq f_1(n) + 2\delta.$$
Using Lemmas \ref{l:prod0} and \ref{l:prod1}, the probability that conditions (1),(2), (3) hold tends to $1$, hence we have 
$$P^n( (gx, g^{-1}x)_x \leq f(n)) \to 1$$
as $n \to \infty$.
\end{proof}

Finally, we put together the previous estimates and use Lemma \ref{l:tau-formula} to prove that translation length grows linearly and 
loxodromic elements are generic. 

\begin{theorem}[Linear growth of translation length] \label{T:gen-lox}
Let $(G, \Gamma)$ be an almost semisimple, nonelementary graph structure, and $L$ the smallest drift given by eq. \eqref{E:small-drift}.
Then for any $\epsilon > 0$ we have 
\[
\frac{\#\{g \in S_n : \tau_X(g) \ge n(L-\epsilon) \}}{\#S_n} \to 1,
\]
as $n\to \infty$. As a consequence, 
\[
\frac{\#\{g \in S_n : g \; \mathrm{is} \; X - \mathrm{loxodromic} \}}{\#S_n} \to 1,
\]
as $n\to \infty$. 
\end{theorem}

\begin{proof}
If we set $f(n) = \eta n$ with $\eta > 0$, then by Proposition \ref{p:g-prod} and Theorem \ref{th:gen_drift} the events 
$(gx, g^{-1}x)_x \leq \eta n$ and $d(x, gx) \geq n(L-\eta)$ occur with probability ($P^n$) which tends to $1$, hence 
by Lemma \ref{l:tau-formula} 
\begin{align*}
P^n \Big( \tau_X(g) \ge n(L- 3 \eta)  \Big) \ge 
P^n \Big ( d(x, gx) - 2(gx, g^{-1}x)_x + O(\delta) \ge n(L- 3 \eta) \Big) 
\end{align*}
which approaches $1$ as $n \to \infty$. This implies the statement if we choose $\epsilon > 3 \eta$. 
The second statements follows immediately since elements with positive translation length are loxodromic.
\end{proof}

\subsection{Genericity of loxodromics for the Markov chain}

We now remark that a very similar proof yields that loxodromics are generic for $\mathbb{P}$-almost every sample path of the Markov chain. 
More precisely, we have the following (which is a reformulation of Theorem \ref{T:Pgeneric}): 

\begin{theorem} \label{T:markov-lox}
Let $(G, \Gamma)$ be an almost semisimple, nonelementary graph structure, and let $L$ be the smallest drift. 
Then for every $\epsilon > 0$, one has 
\[
\PP \Big(\tau_X(w_n) \ge n (L-\epsilon) \Big) \to 1,
\]
as $n \to \infty$. As a consequence, 
\[
\PP \Big(w_n \textup{ is loxodromic on }X \Big) \to 1
\]
as $n \to \infty$.
\end{theorem}

\begin{proof}
The proof is very similar to the proof of Theorem \ref{T:gen-lox}, so we will just sketch it. First, by using the Markov property we establish 
that
$$\lim_{n \to \infty} \mathbb{P} \big((i_n x, t_n^{-1} x)_x \geq g(n) \big) = 0$$
for any choice of function $g : \mathbb{N} \to \mathbb{R}$ such that $\lim_{n \to +\infty} g(n) = + \infty$. Then, by using positivity of the drift 
as in the proof of Lemma \ref{l:prod1} we prove that for each $\epsilon > 0$, we have 
$$\lim_{n \to \infty} \mathbb{P} \big((w_n^{-1} x, t_n^{-1} x)_x \leq n(L - \epsilon)/2 \big) = 0 $$
and 
$$\lim_{n \to \infty} \mathbb{P} \big( (i_n x, w_n x)_x \leq n (L - \epsilon)/2 \big) = 0$$
From the previous three facts, using Lemma \ref{l:fellow_travel} one proves: 
$$\lim_{n \to \infty} \mathbb{P} \Big( (w_n x, w_n^{-1} x)_x \geq f(n) \Big) = 0$$
for any $f : \mathbb{N} \to \mathbb{R}$ such that $\lim_{n \to +\infty} f(n) = +\infty$. 
The theorem then follows immediately from this fact and Corollary \ref{cor:drift_in_measure}, applying the formula of Lemma \ref{l:tau-formula}.  
\end{proof}

\section{Thick graph structures}

\begin{definition}
A graph structure $(G, \Gamma)$ is \emph{thick} if for every vertex $v$ of maximal growth there exists a finite set $B \subseteq G$ such that 
\begin{equation} \label{eq:thick}
G = B \Gamma_v B
\end{equation}
where $\Gamma_v$ is the loop semigroup of $v$. 

In greater generality, if $H < G$ is a subgroup, we say that the graph structure $(G, \Gamma)$ 
is \emph{thick relatively to }$H$ if for any vertex $v$ of maximal growth there exists a finite set $B \subseteq G$ such that
\begin{equation} \label{eq:thick}
H \subseteq B \Gamma_v B.
\end{equation}
\end{definition}

\subsection{The case of only one non-trivial component}

We say that a component $C$ is \emph{non-trivial} if there is at least one closed path of positive length entirely contained in $C$.

\begin{proposition} \label{P:rec-thick}
If a graph structure $(G, \Gamma)$ has only one non-trivial component, then it is thick.
\end{proposition}

\begin{proof}
Let $C$ be the unique maximal component of $\Gamma$.
Every finite path $\gamma$ in the graph can be written as $\gamma = h_1 g h_2$, where $h_1$ is a path from the initial vertex to $C$, $g$ is a path entirely in $C$, and 
$h_1$ is a path going out of $C$. By construction, the lengths of $h_1$ and $h_2$ are uniformly bounded. Fix some vertex $v$ of $C$ and let $s$ be a shortest path from $v$ to the last vertex of $h$. Further, let $t$ be a shortest path from the last vertex of $g$ to $v$.
Then one can write 
$$\gamma = h_1gh_2 = h_1s^{-1} (s g t) t^{-1} h_2$$ 
where $h_{1}s^{-1}$ and $t^{-1} h_2$ vary in a finite set, and $sg t \in \Gamma_v$. Hence $G = B \Gamma_v B$ with $B$ a finite set.
\end{proof}

\subsection{Thick implies nonelementary}

\begin{proposition} \label{P:T->NE}
Fix an action $G \curvearrowright X$ of $G$ on a hyperbolic metric space $X$.
Let $(G, \Gamma)$ be an almost semisimple graph structure, and $H < G$ a nonelementary subgroup. 
If $(G, \Gamma)$ is thick relatively to $H$, then it is nonelementary, i.e.
for any maximal vertex $v$ the action of the loop semigroup $\Gamma_v$ on $X$ is nonelementary.
\end{proposition}

\begin{proof}
Since the action of $H$ is nonelementary, 
there exists a free subgroup $F \subseteq H$ of rank $2$ which quasi-isometrically embeds in $X$. Hence, the orbit map $F \to X$ extends to an embedding $\partial F \to \partial X$, and we identify $\partial F$ with its image.
Thickness implies $F \subseteq B \Gamma_v B$, and taking limit sets in $\partial X$ we see that
\[
\partial F \subset \bigcup_{b \in B} b \cdot \Lambda_{\Gamma_v},
\]
from which we conclude that $\Lambda_{\Gamma_v}$ is infinite. To complete the proof that $\Gamma_v$ is nonelementary, it suffices to show that $\Gamma_v$ does not have a fixed point on $\partial X$.  Suppose toward a contradiction that $p \in \partial X$ is such a fixed point.

Let us write $F = \langle f, g \rangle$ where $f,g$ are free generators of $F$, and consider the sequence of elements $h_{i,j} = f^ig^j$ in $F$. For each $i,j$ there are $a_{i,j},c_{i,j} \in B$ such that $h_{i,j} = a_{i,j}l_{i,j}c_{i,j}$ for some $l_{i,j}$ in $\Gamma_v$. Since $B$ is finite, we may pass to a subsequence and assume that $a_{i,j}=a$ and $c_{i,j} = c$ for all $i,j$. Then $l_{i,j} = a^{-1}h_{i,j}c^{-1}$ fixes the point $p$ for all $i$ and so 
\[
h_{i,j}(c^{-1}(p)) = a(p)
\]
for all $i,j$. Hence $h_{i_0,j_0}^{-1}h_{i,j} = g^{-i_0}f^{i-i_0}g^{j}$ is a sequence of elements of $F$ which fix the point $q = c^{-1}(p) \in \partial F \subset \partial X$. Since $F$ is a free group, this implies that $g^{-i_0}f^{i-i_0}g^{j}$ agree up to powers for infinitely many $i,j$, a clear contradiction.
\end{proof}

From Proposition \ref{P:T->NE} and Theorem \ref{T:gen-lox} we get:

\begin{theorem} \label{T:thick-main}
Let $G \curvearrowright X$ be a nonelementary action of a countable group on a separable, hyperbolic metric space. Suppose that $G$ 
has an almost semisimple graph structure $\Gamma$ which is thick with respect to a nonelementary subgroup $H$. Then loxodromic elements are generic:  
$$\lim_{n \to \infty} \frac{\# \{g \in S_n \ : \ g \textup{ is loxodromic on }X \} }{\#S_n} = 1$$
In fact, the translation length generically grows linearly: there exists $L > 0$ such that 
$$\lim_{n \to \infty} \frac{\# \{g \in S_n \ : \ \tau_X(g) \geq L n  \} }{\#S_n} = 1.$$
\end{theorem}

This completes the proof of Theorem \ref{T:main} in the introduction.

\section{Relative growth quasitightness} \label{sec:rel_quasi}
Fix a graph structure $(G,\Gamma)$. In practice, we will often show that the graph structure  is thick by establishing the property of growth quasitightness. This property was introduced in \cite{arzhantseva2002growth} and and further studied in \cite{Wenyuan-scc}. Our notion of quasitightness depends on the particular graph structure.

Given a path $\gamma$ in $\Gamma$, we say it $c$--almost contains an element $w\in G$ if $\gamma$ contains a subpath $p$ such that 
$w = a \cdot \ev(p) \cdot b$ in $G$, with $|a|, |b| \leq c$. 
We denote as $Y_{w, c}$ the set of paths in $\Gamma$ starting at the initial vertex which do not $c$--almost contain $w$.

\begin{definition}
A graph structure $(G, \Gamma)$ is called \emph{growth quasitight} if there exists $c > 0$ such that for every $w \in G$ the set $Y_{w, c}$ has density zero with respect to $P^n$; that is,
\[
P^n(Y_{w,c}) \to 0 \quad \text{ as } n \to 0.
\]

More in general, given a subgroup $H < G$ we say that $(G, \Gamma)$ is \emph{growth quasitight relatively to }$H$ if there exists a constant $c > 0$ 
such that for every $w \in H$ the set $Y_{w, c}$ has density zero.
\end{definition}

\subsection{Growth quasitight implies thick}

\begin{proposition} \label{L:Hmax}
Let $(G, \Gamma)$ be an almost semisimple graph structure, and $H < G$ a subgroup. 
If $(G, \Gamma)$ is growth quasitight relatively to $H$, then it is thick relatively to $H$.
\end{proposition}

\begin{proof}
Let $C$ be a component of maximal growth, let $v$ a vertex in $C$, and let $\gamma$ be some path from the initial vertex to $v$. Denote the length of $\gamma$ by $d$. 
Let $w \in H$. By growth quasitightness plus maximal growth, there is a path of the form $\gamma\gamma_1$, which $c$--almost contains $w$ and where $\gamma_1$ is entirely contained in $C$. 
Since $\gamma$ has length $d$, the path $\gamma_1$ $(c+d)$--almost contains $w$; that is, 
$$\gamma_1 = p_1p_2p_3$$ 
where $\ev(p_2) = awb$ for $|a|,|b| \le c+d$. Let $q_1$ be a shortest path from $v$ to the initial vertex of $p_2$ and $q_2$ be a shortest path from the terminal vertex of $p_2$ to $v$. Then 
\[
w = (a^{-1} \ev(q_1)^{-1}) \cdot \ev(q_1p_2q_2) \cdot (\ev(q_2)^{-1}b^{-1}),
\]
where $(a^{-1} \ev(q_1)^{-1})$ and $(\ev(q_2)^{-1}b^{-1})$ vary in a finite set $B$.  Since $\ev(q_1p_2q_2) \in \Gamma_v$, this completes the proof.
\end{proof}

Combining Proposition \ref{L:Hmax} with Theorem \ref{T:thick-main} we get: 

\begin{theorem} \label{T:QT-main}
Let $(G, \Gamma)$ be an almost semisimple graph structure which is growth quasitight with respect to a nonelementary subgroup $H$. 
Then loxodromic elements are generic. 
\end{theorem}

This completes the proof of Theorem \ref{T:main}. 

\section{Infinite index subgroups have zero density}

In this section, we prove that in our general setup a subgroup $H < G$ of infinite index has zero density with respect to counting. 
Combined with what we are going to prove in sections \ref{S:rel-hyp} and \ref{RAAGs}, this immediately implies Theorem \ref{sample2} in the introduction. 
Recall that $\ev : \Omega_0 \to G$ is the evaluation map for paths starting at $v_0$. 

\begin{theorem} \label{infiniteindex}
Let $(G,\Gamma)$ be an injective, almost semisimple, thick graph structure.
Let $H<G$ be an infinite index subgroup. Then
\[
P^n \left(\{p\in \Omega_0 : \ev(p) \in H \}\right) \to 0,
\]
as $n\to \infty$. That is, the proportion of paths starting at $v_0$ which spell an element of $H$ goes to $0$ as the length of the path goes to $\infty$.
\end{theorem}

The proof is an adaptation of (\cite{gouezel2015entropy}, Theorem 4.3) to the non-hyperbolic case.
We will consider an extension $\Gamma_{H}$ of $\Gamma=(V,E)$ defined as follows. The vertex set of $\Gamma_{H}$ is $V\times H \setminus G$.
For any edge $\sigma:x\to y$ in $\Gamma$ 
there is an edge in $\Gamma_{H}$ from $(x,Hg)$ to $(y,Hg g')$ where $g' =  \textup{ev}(\sigma)$.

\begin{lemma}\label{infdiam}
Let $C$ be a component of maximal growth. For any $v_{1}\in C$ and $g_{1}\in G$ there are infinitely many $Hg\in H\setminus G$ such that $(v_{1},Hg)$ can be reached from $(v_1,Hg_{1})$ by a  path contained in $C\times H\setminus G$.
\end{lemma}

\begin{proof}
Suppose not, so that the only points of $H\setminus G$ that can be reached in this manner are 
$\{ Hz:z\in T \}$ where $T$ is a set of size $D$. Consider $w\in G$. 
By thickness, there exists a finite set $B \subseteq G$ and some path 
$\gamma$ lying in $C$, starting and ending at $v_{1}$ such that 
$$\ev(\gamma)=g_{2}wg_{3}$$ 
where $g_{2},g_{3}$ lie in $B$. 
Then $\gamma$ lifts to a path in $\Gamma_{H}$ from $(v_{1},Hg_{1})$ to $(v_{1}, Hg_{1}\textup{ev}(\gamma))$.
By assumption, this implies $Hg_{1} \textup{ev}(\gamma)=Hz$ for some $z\in T$. Thus, there is an $h\in H$ with 
$g_{2}wg_{3}=\textup{ev}(\gamma)= g^{-1}_{1} h z$ and hence $w\subset B^{-1}g^{-1}_{1}HTB^{-1}$.
Thus there is a finite subset $\Upsilon=B^{-1} g^{-1}_{1}\cup TB^{-1}$ with $G=\Upsilon H \Upsilon$, so by Neumann's theorem \cite{Neumann} $H$ must be of finite index, giving a contradiction.
\end{proof}

The following general result about Markov chains is Lemma 4.4 of \cite{gouezel2015entropy}. 
\begin{lemma}\label{Markov}
Let $X_n$ be a Markov chain on a countable set $V$, and $m$ a stationary measure. Let $\tilde{V}$ be the set of points $x\in V$ such that $\sum_{y:x\to y} m(y)=\infty$ where $x\to y$ means there is a positive probability path from $x$ to $y$. Then for all $x\in V$ and $x'\in \tilde{V}$ we have $P_{x}(X_{n}=x')\to 0$.
\end{lemma}

Combining Lemmas \ref{Markov} and \ref{infdiam} we obtain:
\begin{corollary}\label{pointwisedecay}
For any $x_{1},x_{2}\in \Gamma$ lying in a maximal component $C$ and $g_{1},g_{2}\in G$, 
the number of paths of length $n$ in $\Gamma_H$ from $(x_{1},Hg_{1})$ to $(x_{2},Hg_{2})$ is $o(\lambda^{n})$. 
\end{corollary}

\begin{proof}
The Markov chain $\mu$ on $\Gamma$ restricts to a Markov chain $\mu_{C}$ on $C$, which in turn lifts to a Markov chain $\mu_{C,H}$ on the induced graph $C_{H}$ on the vertex set $C\times H\setminus G \subset \Gamma_{H}$ (obtained by assigning to an edge the transition probability of its projection to $C$.
A $\mu_{C,H}$ stationary measure $\tilde{m}$ on $C_H$ is given by taking the product of the  stationary measure $m$ on $C$ and the counting measure on $H\setminus \Gamma$. 
Any vertex $v\in C$ has positive $m$ measure and all lifts of $v$ in $C_{H}$ have equal positive $\tilde{m}$ measure. Thus, Lemma \ref{infdiam} implies 
$\sum_{y:x\to y}\tilde{m}(y)\to \infty$.
The corollary now follows by applying Lemma \ref{Markov} to the chain $\mu_{C,H}$. 
\end{proof}
Note, paths of length $n$ in $\Gamma_H$ from $(x_{1},Hg_{1})$ to $(x_{2},Hg_{2})$ are in bijection with paths of length $n$ in $\Gamma$ beginning at $x_1$, ending at $x_2$, and evaluating to elements of $g^{-1}_{1}H g_{2}$. Thus, we obtain:
\begin{corollary}\label{smallcount}
For any $x_{1},x_{2}\in \Gamma$ lying in a maximal component $C$ and $g_{1},g_{2}\in G$, the number of paths of length $n$ in $\Gamma$ beginning at $x_1$, ending at $x_2$, and evaluating to elements of $g^{-1}_{1}H g_{2}$ is $o(\lambda^n)$.
\end{corollary}

We now complete the proof of Theorem \ref{infiniteindex}. Given $k>0$, let $P_{n,k}$ (resp. $Q_{n,k}$) be the set of paths $p \in \Omega_0$ of length $n$ which spend time at most $k$ (resp. more than $k$) in non-maximal components.

Note that there is a $\eta<\lambda$ with $|Q_{n,k}|\leq \eta^{k}\lambda^{n-k}$ for all $n$ and $k$.
Now, consider a path $\gamma$ in $P_{n,k}\cap \ev^{-1}H$. We can decompose it as $\gamma=\gamma_{1}\gamma_{2}\gamma_{3}$ where $\gamma_{1}$ and $\gamma_{3}$ have length adding up to at most $k$ and $\gamma_{2}$ is contained in a maximal component $C$.
Since a path in $P_{n,k}$ spends at most time $k$ in nonmaximal components, there are only $D^{k}$ possibilities for $\gamma_{1}$ and $\gamma_{3}$, 
where $D$ depends only on the graph.
On the other hand, by Corollary \ref{smallcount}, for a path in $\ev^{-1}H$, there are at most $f(n) = o(\lambda^n)$ possibilities for 
 $\gamma_2$, thus $|P_{n,k}\cap \ev^{-1}H|\leq D^{k}f(n)$ for all $k<n$ where $f(n)=o(\lambda^n)$.
Hence, 
$$P^{n}(\ev^{-1}H)\leq C''\lambda^{-n}(|P_{n,k}\cap \ev^{-1}H|+|Q_{n,k}|)\leq C''D^{k}\lambda^{-n}f(n)+C''(\eta/\lambda)^{k}.$$
Fixing $k$ we see that 
$$\lim \sup_{n\to \infty}P^{n}(\ev^{-1}H)\leq C''(\eta/\lambda)^{k}.$$
As this is true for arbitrary $k$, we get $\lim_{n \to \infty} P^{n}(\ev^{-1}H) = 0$, as claimed.

\section{Application to relatively hyperbolic groups} \label{S:rel-hyp}

In this section, we show how our main theorem applies to a large class of relatively hyperbolic groups. 

Let $G$ be a finitely generated group, and $\mathcal{P}$ be a collection of subgroups. Following \cite{bowditch2012relatively}, let us recall that  
$G$ is \emph{hyperbolic relative to} $\mathcal{P}$ if there is a compactum $M$ on which $G$ acts geometrically finitely, and the 
maximal parabolic subgroups are the elements of $\mathcal{P}$. Such a compactum $M$ is then unique up to $G$-equivariant 
homeomorphisms, and it is called the \emph{Bowditch boundary} of $G$. We will denote it as $\partial G$.

\medskip
More precisely, let $G$ act by homeomorphisms on a compact, perfect, metrizable space $M$. Then a point $\zeta \in M$ is called \emph{conical} if there is a sequence $(g_n)$ and distinct points $\alpha, \beta \in M$ such that 
$g_n \zeta \to \alpha$ and $g_n \eta \to \beta$ for all $\eta \in M \setminus \{ \zeta \}$. A point $\zeta \in M$ is called \emph{bounded parabolic} if 
the stabilizer of $\zeta$ in $G$ is infinite, and acts cocompactly on $M \setminus \{ \zeta \}$.
We say that the action of $G$ on $M$ is a \emph{convergence action} if $G$ acts properly discontinuously on triples of elements of $M$, and the action is \emph{geometrically finite} if it is a convergence action and every point of $M$ is either a conical limit point or a bounded parabolic point. 
Note that there are only countably many parabolic points. Finally, the \emph{maximal parabolic subgroups} are 
the stabilizers of bounded parabolic points. 
We refer the reader to \cite{farb1998relatively, bowditch2012relatively, osin2006relatively, groves2008dehn} for the relevant background material.

Fix a relatively hyperbolic group $(G, \P)$, a generating set $S$, and let $d_G$ denote distance in $G$ with respect to $S$. 
Let $\hat G$ be the vertices of $\cay(G,S \cup \P)$ with the induced metric, which we denote by $\hat d$ or $d_{\hat G}$. Here $\cay(G, S \cup \P)$ is the corresponding electrified Cayley graph, that is the Cayley graph of $G$ with respect to the generating set $S \cup \bigcup \P$. We remind the reader that $\cay(G, S \cup \P)$ is hyperbolic and that $\partial \cay(G,S \cup \P)$ naturally includes as a subspace into $\partial G$, the complement of which is the collection of parabolic fixed points. 
%\marpar{Do these facts appear in Bowditch?}

Following \cite{YangPS}, the Bowditch boundary $\partial G$ is equipped with a quasiconformal, nonatomic measure $\nu$, which is given by the Patterson-Sullivan 
construction by taking average on balls \emph{for the word metric} on $G$.

\begin{definition}
We define a relatively hyperbolic group $G$ to be \emph{pleasant} if its action on $\partial G \times \partial G$ with the measure $\nu \times \nu$ 
is ergodic. 
\end{definition}

We will also see (Proposition \ref{doublyergodic}) that a relatively hyperbolic group is pleasant if it admits a geometrically finitely action on a $CAT(-1)$ proper metric space. 
For instance, geometrically finite Kleinian groups satisfy this hypothesis. Note that, once $G$ admits such an action, the theorem works for isometric actions of $G$ on \emph{any} hyperbolic, metric space $X$. 

In this section, we will prove the following result. 

\begin{theorem}
Let $G$ be a pleasant, relatively hyperbolic group, and let $(G, \Gamma)$ be a geodesic combing. Then 
for any nonelementary action of $G$ on a hyperbolic metric space $X$, the graph structure $(G, \Gamma)$ is 
nonelementary. 
\end{theorem}

Combining this result with Theorem \ref{T:main}, the discussion in section \ref{S:introRH}, and the fact that relatively hyperbolic groups 
has pure exponential growth for any generating set (\cite{YangPS}, Theorem 1.9) this establishes Theorem \ref{T:mainRH} in the introduction. 

In fact, using very recent work of W. Yang \cite{Wenyuan-scc}, the theorem may be extended to all nontrivial relatively hyperbolic groups, as relatively hyperbolic groups contain strongly contracting elements by \cite{Arzhantseva-Cashen-Tao} and \cite{Wenyuan-scc} (see Corollary \ref{C:wenyuan} in the introduction). However, we give a self-contained argument here. 

\subsection{Fellow traveling in the Cayley graph and coned-off space}
We will need the following proposition, which is certainly known to experts. We provide a proof for completeness.

\begin{proposition} \label{prop:stable_overlap}
For $K,C \ge0$, there are $D,L \ge0$ such that the following holds. Suppose that $\gamma = [a,b]$ is a geodesic in $\cay(G,S)$ with length at least $L$ which projects to a $K$--quasigeodesic in $\cay(G,S \cup \P)$. Let $\gamma'$ be any other geodesic in $\cay(G,S)$ whose endpoints have distance no more than $C$ from $a,b$ in $\cay(G, S \cup \P)$.
Then there are $a',b' \in \gamma'$ such that 
\[
d_G(a,a') \le D \quad \text{and} \quad d_G(b,b') \le D.
\]
\end{proposition}

We will use the following theorem of Osin:

\begin{theorem}[\cite{osin2006relatively}, Theorem 3.26] \label{th:Osin_thin}
There is an $\nu\ge0$ such that if $p, q, r$ are sides of a geodesic triangle in $\cay(G, S \cup \P)$, then for any vertex $v$ on $p$ there exists a vertex $u$ on either $q$ or $r$ such that
\[
d_G(v,u) \le \nu.
\]
\end{theorem}

\begin{proof}[Proof of Proposition \ref{prop:stable_overlap}]
Suppose that $d_G(a,b) > K(2C+4\nu+K)$ so that $d_{\hat G}(a,b) > 2C+4\nu$, for $\nu$ as in Theorem \ref{th:Osin_thin}. 

Let $c,c'$ be geodesics in $\cay(G,S\cup \P)$ joining the endpoints of $\gamma, \gamma'$ respectively. Note that by assumption the initial and terminal endpoints of these geodesics are at $\hat d$--distance less than $C$ from one another.  Pick vertices $c_a,c_b$ on $c$ (ordered $a,c_a,c_b,b$) so that $d_{\hat G}(a,c_a)= d_{\hat G}(b,c_b)= C+2\nu$. (This is possible since $d_{\hat G}(a,b) > 2C+4\nu$.)
Consider a geodesic quadrilateral with opposite sides $c, c'$. Applying Theorem \ref{th:Osin_thin} twice, we may find vertices $c'_a,c'_b \in c'$ such that $d_G(c_a,c_a')\le 2\nu$ and  $d_G(c_b,c_b') \le 2\nu$. 

Now using, for example, (\cite{hruska}, Lemma 8.8), we can
 find vertices $\gamma_a,\gamma_b \in \gamma$ and $\gamma'_a,\gamma'_b \in \gamma'$ which have $d_G$--distance at most $L$ from $c_a,c_b,c'_a,c'_b$, respectively, where $L\ge 0$ depends only on $(G,\P)$ and $S$. Note that $d_G(\gamma_a,\gamma'_a)\le 2(L+\nu)$ and $d_G(\gamma_b,\gamma'_b) \le 2(L+\nu)$. Moreover, 
\[
d_{\hat G}(a,\gamma_a) \le d_{\hat G}(a,c_a) + d_{\hat G}(c_a,\gamma_a) \le C+2\nu+L,
\]
and so $d_G(a, \gamma_a) \le K(C+2\nu+L+1)$ since both $a$ and $\gamma_a$ occur along $\gamma$. Similarly, $d_G(b,\gamma_b) \le K(C+2\nu+L+1)$.

Putting everything together, after setting $a' = \gamma'_a$ and $b' = \gamma'_b$, we see that each of $d_G(a,a')$ and $d_G(b,b')$ are less than  $2(L+\nu)+ K(C+2\nu+L+1)$ and this completes the proof.
\end{proof}

\subsection{Patterson-Sullivan measures and sphere averages} \label{S:PS-RH}

Continuing with the notation from the previous section,
let $h$ be the exponent of convergence for $\cay(G,S)$. That is,
\[
h = \lim _{n \to \infty} \log |B_n| /n,
\]
where $B_n$ denotes the ball of radius $n$ in $G$ with respect to $d_G$.

\begin{definition} 
For $g\in G$, we define the \emph{large shadow} $\Pi_{r}(g)$ at $g$ to be the set of $\zeta \in \partial G$ such that  there exists \emph{some} geodesic in $\cay(G,S)$ from $1$ converging to $\zeta$ intersecting $B_{r}(g)$.
Similarly, the \emph{small shadow} $\pi_{r}(g)$ is the set of $\zeta \in \partial G$ such that \emph{every} geodesic in $\cay(G,S)$ from $1$ converging to $\zeta$ intersects $B_{r}(g)$.
\end{definition}

In Theorem 1.7 and Proposition A4 of \cite{YangPS} Yang constructs an $h$-quasiconformal ergodic density $\nu$ without atoms for the word metric on the Bowditch boundary $\partial G$. 
In Lemma 4.3 of \cite{YangPS} he shows that this satisfies the shadow lemma: for large enough $r$:
\begin{equation} \label{E:shadowlemma} 
\nu(\pi_{r}(g)) \simeq \nu(\Pi_{r}(g))\simeq e^{-h \cdot d_G(1,g)}
\end{equation}
(up to a uniform multiplicative constant). In particular, $\nu$ has full support on $\partial G$. In what follows, $S_n$ denotes the set of elements $g \in G$ with $d_S(1,g)=n$.

\begin{lemma} \label{L:boundarymeasure}
There is a $C>0$ such that for any Borel set $A\subset G\cup \partial G$ one has
$$\lim \sup_{n\to \infty}\frac{|A\cap S_{n}|}{|S_{n}|}\leq C\nu(\overline{A}),$$
where $\overline A$ denote the closure of $A$ in $G \cup \partial G$.
\end{lemma}

\begin{proof}
Let $A\subset G\cup \partial G$ be a Borel set. Since the number of elements in a ball of radius $r$ in $\cay(G,S)$ is universally bounded, a point of $\partial G$ lies in at most $D$ small shadows $\pi_{r}(g), g\in S_{n}$ where $D$ depends only on $r$. Thus,
$$\sum_{g\in S_{n} \cap A}\nu\left(\pi_{r}(g)\right)\leq D \nu\left(\bigcup_{g\in S_{n}\cap A}\pi_{r}(g)\right)$$
Moreover, if we denote $A_n := A \setminus B_{n-1}$ then $A_{n+1} \subseteq A_n$ and 
$$\bigcap_{n\in \N}\bigcup_{g\in A_n}\pi_{r}(g)\subset \overline{A}.$$
Indeed, if $\zeta \in \bigcap_{n\in \N}\bigcup_{g\in A_n}\pi_{r}(g)$, then there are $g_n \in A$ with $|g_n| \ge n$ such that some (any) geodesic from the identity to $\zeta$ meets $B_r(g_n)$. Hence, $g_n \to \zeta$ and so $\zeta \in \overline A$.

Thus, since $S_n \cap A \subseteq A_n$ we have for large enough $n$
$$\nu\left(\bigcup_{g\in S_{n}\cap A}\pi_{r}(g) \right)\leq \nu\left(\bigcup_{g\in A_n}\pi_{r}(g)\right) \leq 2\nu(\overline{A})$$
so by exponential growth and the shadow lemma \eqref{E:shadowlemma}
$$\frac{|S_{n}\cap A|}{|S_{n}|}\simeq e^{-hn}|S_{n}\cap A|\simeq \sum_{g\in S_{n}\cap A}\nu(\pi_{r}(g))  \lesssim 2 D \nu(\overline{A}).$$
\end{proof}

\subsection{Growth quasitightness for relatively hyperbolic groups}
We will now establish a form of relative growth quasitightness for a relatively hyperbolic group $G$. 

Let $w$ be an element of $G$. A \emph{$w$-path} is an infinite path of the form $l_w = \bigcup_{i \in \Z} w^i \gamma_w$ in the Cayley graph $\cay(G, S)$, where $\gamma_w = [1, w]$ is a geodesic segment joining the identity and $w$. Of course there may be finitely many choices of $l_w$ for each $w$.

\begin{definition}
The element $w$ is called \define{$K$--bounded} if some $w$-path $l_w$ (with the arc length parameterization) in the Cayley graph $\cay(G,S)$ projects to a $K$-quasigeodesic in the electrified graph $\cay(G, S \cup \P)$.
\end{definition}

The following lemma is well-known. See for example \cite{DMS, ADT}.

\begin{lemma} \label{lem:stable}
For each $K$, there is a function $f \colon \N \to \N$ such that if $w$ is $K$--bounded, then every $w$-path $l_w$  
is an $f$-stable quasigeodesic in the Cayley graph $\cay(G,S)$. 
\end{lemma}
Recall that $l_w$ being $f$--stable means that any $K$--quasigeodesic with endpoints on $l_w$ has Hausdorff distance at most $f(K)$ from the subpath of $l_w$ its endpoints span.

Given $w \in G$ and $c \geq 0$, we say that a (finite or infinite) geodesic  $\gamma$ \emph{$c$-almost contains} $w$ if there exists $g \in G$ such that $d_G(g, \gamma) \leq c$ and 
$d_G(gw, \gamma) \leq c$. 
Let $X_{w,c}$ be the set of $h\in G$ such that there exists a geodesic $\gamma$ from identity to $h$ which does not $c$-almost contain $w$. 
That is, for every $g\in N_c(\gamma)$, $\gamma$ does NOT pass within distance $c$ of $gw$.

\begin{proposition} \label{P:gqt}
For each $K\ge 1$, there is $c\ge0$ such that for every $K$--bounded $w \in G$ we have 
$$\frac{|B_n \cap X_{w,c}|}{e^{h n}}\to 0$$ 
as $n \to \infty$. 
\end{proposition}

We remark that, for fixed $c$, it suffices to prove the proposition for sufficiently long $w$, that is, where $|w|_S$ is sufficiently large.
We will prove this proposition by using the ergodicity of the double boundary (Proposition \ref{doublyergodic}).
To do this, we will apply Proposition \ref{prop:stable_overlap} several times, for $K$ the boundedness constant. Hence, we fix $K$ once and for all, and consider the constant $D$ produced by that proposition as a function of $C$ alone and write $D = D(C)$.

\medskip
Let $Z_{w,c}$ be the set of pairs $(\alpha,\beta)$ in $\partial G \times \partial G$ such that for every bi-infinite geodesic $\gamma$ in $\cay(G,S)$ joining $\alpha$ and $\beta$ there exist infinitely many $x \in N_c(\gamma)$ such that $\gamma$ passes within $c$ of $xw$.
Let $Z^{n}_{w,c}$ be the set of pairs $(\alpha,\beta)$ in $\partial G \times \partial G$ such that for every bi-infinite geodesic $\gamma$ joining $\alpha$ and $\beta$  there are at least  $n$ 
elements $x \in N_c(\gamma)$ such that $\gamma$ passes within $c$ of $xw$.
By definition,
$Z_{w,c}=\bigcap_{n\in \N}Z^n_{w,c}$.
Moreover, for each $n,w,c$, the sets $Z_{w,c}$ and $Z^{n}_{w,c}$ are $G$-invariant subsets of  $\partial G \times \partial G$.

Furthermore we have
\begin{lemma} \label{prop:biinfinite}
For each $K\ge0 $ there is a constant $c_0=c_0(K)$ such that for all $c \ge c_0$, $Z_{w,c}$ contains a pair of conical points
for every $K$--bounded $w \in G$.
\end{lemma}

\begin{proof}
Let $f$ be the function given by Lemma \ref{lem:stable}. By definition, the $w$-path $l_w$ projects to a $K$-quasigeodesic in $\cay(G, S \cup \P)$, 
hence it has two distinct limit points $(w^{-\infty},w^\infty)$  in the Bowditch boundary $\partial G$. Then, by connecting further and further points on $l_w$ by a geodesic in the $\cay(G, S)$, using $f$-stability and the Ascoli-Arzel\'a theorem, one constructs a geodesic in $\cay(G, S)$ which connects $w^\infty$ to $w^{-\infty}$ and $c$-fellow travels $l_w$, once $c\ge 2f(1)$. Hence, $(w^{-\infty},w^\infty) \in Z_{w,c}$.
\end{proof}

\begin{lemma} \label{prop:nonempy_int}
For each $K$, there is a $c_1 = c_1(K)$ and $L_1 = L_1(K)$ such that for $c\ge c_1$ and for any $K$--bounded $w\in G$ with $|w| \geq L_1$, the set $Z^{n}_{w,c}$ has nonempty interior. More precisely, the interior of
$Z^{n}_{w,c}$ contains every pair of conical points in $Z^{n}_{w,c_0}$ where $c_0$ is as in Lemma \ref{prop:biinfinite}.
\end{lemma}

\begin{proof}
Suppose that $(\alpha,\beta) \in Z^{n}_{w,c_0}$ is a pair of conical points, and pick a geodesic $\gamma$ joining $\alpha$ and $\beta$. Then by definition there are segments $[x_j, x'_j] \subseteq [\alpha, \beta]$ for  $1\le j\le n$ and points $z_j$ such that 
$d_G(x_j,z_j), d_G(x'_j,z_jw) \le c_0$ (for $1\le j\le n$). 

Now, let $\alpha_i, \beta_i \in \partial G$ such that $\alpha_i \to \alpha$ and $\beta_i \to \beta$. 
In particular, for some uniform $R$, the projections of the geodesics $[\alpha_i, \beta_i]$ to $\cay(G,S \cup \P)$ $R$--fellow travel the projection of $[\alpha,\beta]$ for longer and longer intervals. Hence, there is an $N\ge0$ so that for $i\ge N$, each $[\alpha_i, \beta_i]$ passes with $\hat d$--distance $R$ from $[x_j, x_j'] \subset [\alpha,\beta]$.

Since $w$ is $K$--bounded, the geodesic segment $z_j \cdot [1,w] = [z_j, z_j w]$ projects to a $K$--quasigeodesic, so we may apply Proposition \ref{prop:stable_overlap} 
(with $C = R+c_0$) to find constants $D$, $L$ such that if $|w| \geq L$ there exist points $y_j,y_j' \in [\alpha_i,\beta_i]$ such that $d_G(y_j,z_j), d_G(y_j',z_jw) \le D$.

Setting $c_1 = D$, $L_1 = L$ we see that for sufficiently large $i\ge0$, $(\alpha_i,\beta_i) \in Z^{n}_{w,c_1}$ and this completes the proof.
\end{proof}

Since $G$ is pleasant, the action of $G$ on $\partial G \times \partial G$ is ergodic, hence Lemma \ref{prop:nonempy_int} implies

\begin{lemma}\label{closedmeasure} \label{intmeasure}
For $c \geq c_1(K)$ and for any $K$--bounded $w\in G$ with $|w| \geq L_1(K)$, the set $Z^{n}_{w,c}$ has full $\nu \times \nu$ measure. Hence, under the same hypotheses the set $Z_{w,c}$ has full $\nu \times \nu$ measure. 
\end{lemma}

\begin{proof}
For each $n$, the set $Z^n_{w, c}$ is $G$-invariant, hence by ergodicity its measure is either $0$ or $1$. Since it has nonempty interior and the measure $\nu \times \nu$ has full support, then it must have full measure. The second claim follows since $Z_{w, c} = \cap_n Z^n_{w, c}$.
\end{proof}

Let $\Lambda_{w,c}\subset \partial G$ be the set of conical points $\alpha \in \partial G$ such that for every geodesic ray $\gamma$ from the identity converging to $\alpha$ there are infinitely many points $g\in N_c(\gamma)$ such that $\gamma$ passes within distance $c$ of $gw$.
Using Proposition \ref{prop:stable_overlap} just as in Lemma \ref{prop:nonempy_int} we have:

\begin{lemma}
For each $K\ge 0$ there is a $c_2 = c_2(K)$ and $L_2 = L_2(K)$ such that for $c \ge c_2$, and for any $K$--bounded $w \in G$ with $|w| \geq L_2$, if 
$(\alpha, \beta)\in Z_{w,c_1}$, then either $\alpha$ or $\beta$ is in $\Lambda_{w,c}$. 
\end{lemma}

This implies

\begin{corollary} \label{unilateralmeasurezero}
For each $c \ge c_2$ and for any $K$-bounded $w \in G$ with $|w| \geq L_2$, the set $\Lambda_{w,c}$ has full $\nu$ measure. 
\end{corollary}

\begin{lemma} \label{P:l1}
For each $K\ge 0$ there is a $c_3 = c_3(K)$ and $L_3 = L_3(K)$ such that for each $c\ge c_3$ and for any $K$--bounded $w\in G$ with $|w| \geq L_3$,
the closure of $X_{w,c}$ is contained in $\partial G\setminus \Lambda_{w,c_2}$.
\end{lemma}

\begin{proof}
If this were false, then for all large $c\ge0$ there would be a sequence $(y_i) \subseteq X_{w,c}$ converging to $\eta \in  \Lambda_{w,c_2}$. 
Since $\eta$ is not a parabolic fixed point, 
then one can view $\eta$ as belonging to the boundary of $\cay(G,S\cup \P)$. Then the projections to $\cay (G,S \cup \P)$ of any geodesics $\gamma_i = [1,y_i]$ must $R$--fellow travel the projection to $\cay(G, S \cup \P)$ of $[1,\eta]$ for longer and longer intervals, where $R$ is independent of $w$. If $\eta$ were in $\Lambda_{w,c_2}$, then just as in the proof of Lemma \ref{prop:nonempy_int}, we would obtain by applying Proposition \ref{prop:stable_overlap} (with $C = c_2+R$) two constants $L_3, c_3$ such that for $|w| \geq L$ 
and for large $i$, the geodesic $\gamma_i$ $c_3$--almost contains $w$. 
Hence, for $c\ge c_3$ we obtain a contradiction to $y_i \in X_{w,c}$ for all $i$. This completes the proof.
\end{proof}

We are now in position to prove Proposition \ref{P:gqt}. 

\begin{proof}[Proof of Proposition \ref{P:gqt}]
By Lemma \ref{P:l1} and Corollary \ref{unilateralmeasurezero}, for $c \ge c_3$
$$\nu(\overline{X_{w, c}}) \leq \nu(\partial G \setminus \Lambda_{w, c_2}) = 0$$
Hence, applying Lemma \ref{L:boundarymeasure}, for large enough $c>0$ we have
$$\lim \sup_{n}e^{-h n}|B_n \cap X_{w,c}|\leq C\nu(\overline{X_{w,c}})=0.$$
\end{proof}

\subsection{The loop semigroup is nonelementary}

We will now assume that $(G,\P)$ is a pleasant relatively hyperbolic group which admits a geodesic combing for the generating set $S$.

Recall that by work of Antol\'in-Ciobanu \cite{AntolinCiobanu}, if the parabolic subgroups are geodesically completable, then every generating 
set for $G$ can be extended to a generating set for which $G$ has a geodesic combing. From here on, we will use such a generating set. 

Then for each $w \in G$ and constant $c$, let us recall that $Y_{w, c}$ is the set of paths $\gamma$ in the directing graph from the initial vertex which do not $c$-almost contain $w$, i.e. such that 
one cannot write $\ev(\gamma) = a_1 w a_2$ in $G$, with $|a_i| \leq c$ for $i = 1, 2$. By identifying paths from the identity with group elements, it is immediate from the 
definition that $Y_{w, c} \subseteq X_{w, c}$. Hence, by Proposition \ref{P:gqt}, also $Y_{w, c}$ has zero density if $w$ is $K$--bounded.

\begin{proposition}
Let $(G, \Gamma)$ be a geodesic combing for a pleasant, relatively hyperbolic group. Then $(G, \Gamma)$ is nonelementary. 
\end{proposition}

\begin{proof}
We are going to prove that the graph structure is thick relative to a nonelementary, free subgroup $F < G$, which yields the claim 
by Proposition \ref{P:T->NE}.
Let $v$ be a vertex of maximal growth and $w$ be any $K$--bounded word. Let $d=\textup{diam }\Gamma$. Let $c$ be the constant from Proposition 
\ref{P:gqt}. Let $h_1$ be a group element representing a path from the initial vertex to $v$, and consider the set 
$$\Sigma = \{ h_1h_2 \ : \ h_2 \in \Gamma_v \}$$
Since $v$ has maximal growth and $Y_{w, c}$ has zero density, the set $\Sigma$ contains a path $h$ which does not belong to $Y_{w, c}$. 
Then there is a path $h=h_{1}h_{2}$ such that $h_{1}$ has length $\leq d$, $h_{2}$ is entirely contained in the component $C_{v}$ containing $v$, and $h_{2}$ contains a subpath of the form $w'=awb$ where $a$ and $b$ have length less than $c+d$.
Let $s$ be a path from $v$ to the start of $w'$ and $t$ be a path from the end of $w'$ to $v$, each of length at most $D$.
Then $sw't$ is in $\Gamma_{v}$
and $w=as^{-1}(swt)t^{-1}b\subset B\Gamma_{v}B$ where $B$ is a finite  set.

To complete the proof, it suffices to show that $B\Gamma_vB$ contains a nonelementary subgroup (Proposition \ref{P:T->NE}). Using a standard ping-pong argument, construct a free subgroup $H = \langle f ,g \rangle \le G$ which $K$--quasi-isometrically embeds in $\cay(G,S \cup \P)$ and which $K$--quasi-isometrically embeds into $X$ for some $K\ge0$. (Indeed, by \cite{TT}, a random $2$-generator subgroup of $G$ will have this property.) For this $K$, let $B$ be the finite subset produced above enlarged to contain $f^\pm, g^\pm$. Then for any $w\in H$, at least one of $w$, $wf$, or $wg$ is cyclically reduced in $H$ and hence $K$--bounded in $G$. Hence $w \in B \Gamma_v B$ and so $H \le B \Gamma_v B$, as required.
\end{proof}

\subsection{Double Ergodicity} 

We conclude this section by proving that a group $G$ which admits a geometrically finite action on a CAT$(-1)$ space is pleasant. 

Let us assume that $G$ acts geometrically finitely on a CAT($-1$) space $Y$. Recall that an orbit map $G \to Y$ induces an embedding $\partial G \to \partial Y$ \cite[Theorem 9.4]{bowditch2012relatively}, so we can identify $\partial Y$ with the Bowditch boundary of $G$. We continue to denote the pushforward of the measure $\nu$ to $\partial Y$ by $\nu$.

\begin{proposition}\label{doublyergodic}
Suppose $G$ acts geometrically finitely on a $\mathrm{CAT}(-1)$ space $Y$. Then the action of $G$ on $\partial Y \times \partial Y$ is ergodic with respect to $\nu \times \nu$.
\end{proposition}

We remind the reader that $\nu$ is quasiconformal with respect to the word metric rather than the metric on $Y$.

\begin{proof}
Assume $G$ acts geometrically finitely on a CAT$(-1)$ space $Y$, with elements of $\P$ being the parabolic subgroups. 
Recall that the Bowditch boundary $\partial G$ is identified with the Gromov boundary of Y.
Let $d = d_G$ still denote the word metric on $G$.
For $\zeta \in \partial G$ and $g,h\in G$ let
$$\beta_{\zeta}(g,h)=\lim \sup_{z\to \zeta}\left(d(g,z)-d(h,z)\right)$$ be the Busemann functions for the word metric on $G$.
By W. Yang's Lemma 2.20 in \cite{YangPS} there is a $C>0$ such that
for every conical $\zeta \in \partial G$ we have 
\begin{equation} \label{E:liminf}
|\lim \sup_{z\to \zeta}[d(g,z)-d(h,z)]-\lim \inf_{z\to \zeta}[d(g,z)-d(h,z)]|<C
\end{equation}
Moreover, the Patterson-Sullivan measure $\nu$ (for the word metric) gives full measure to conical points and is a quasiconformal density in the sense that there is a $D>0$ such that:
\begin{equation} \label{E:qc}
D^{-1}e^{-h\beta_{\zeta}(g,e)}\leq \frac{dg_{\star}\nu}{d\nu}(\zeta)\leq D e^{-h\beta_{\zeta}(g,e)}
\end{equation}
for all $g \in G$ and $\zeta \in \partial G$. We claim there is a $G$-invariant measure in the measure class of $\nu \times \nu$.
Indeed, let
$$\rho_{e}(\zeta ,\xi)= \lim \sup_{z\to\zeta, y\to \xi} 
\left( \frac{d(e,y)+d(e,z)-d(y,z)}{2} \right).$$
Define a locally finite measure $m'$ on $(\partial G \times \partial G)\setminus Diag$ by 
$$dm'(\zeta, \xi)=e^{2h \rho_{e}(\zeta,\xi)}\ d\nu(\zeta)\ d\nu(\xi)$$
We now claim that the measure $m'$ is $G$ quasi-invariant with a uniformly bounded derivative.  Indeed, we can compute
$$2 \rho_e(g^{-1} \zeta, g^{-1}\xi) - 2 \rho_e(\zeta, \xi) = $$
$$ = \limsup_{z \to \zeta, y \to \xi} \left[ d(e, g^{-1}y) + d(e, g^{-1}z) - d(g^{-1}y, g^{-1}z) \right] - \limsup_{z \to \zeta, y \to \xi} \left[ d(e, y) + d(e, z) - d(y, z) \right] = $$
$$ = \limsup_{y \to \xi} \left[ d(g, y) - d(e, y) \right] + \limsup_{z \to \zeta} \left[ d(g, z) - d(e, z) \right] + O(1) = \beta_\xi(g, e) + \beta_\zeta(g, e) + O(1)$$
(where we could distribute the limsup since the limsup and liminf are within bounded difference (see \eqref{E:liminf}).
Hence, combining this with \eqref{E:qc} one gets that the Radon-Nykodym cocycle is uniformly bounded, i.e.
$$\frac{dg_\star m'}{dm'} (\zeta, \xi)=e^{2h \rho_{e}(g^{-1} \zeta, g^{-1} \xi)- 2h \rho_{e}(\zeta, \xi)}\ \frac{dg_\star\nu}{d\nu}(\zeta)\ \frac{dg_\star\nu}{d\nu}(\xi) \cong 1$$
Hence, by a general fact in ergodic theory the Radon-Nykodym cocycle is also a coboundary (see \cite{Furman}, Proposition 1). Thus, there exists a $G$-invariant measure $m$ on $(\partial G \times \partial G)\setminus Diag$ in the same measure class as $m'$, hence also in the same measure class as $\nu \times \nu$. By \cite{YangPS}, the Patterson-Sullivan measure is supported on conical limit points. Thus, $m$ is also supported on pairs of conical limit points.
By Theorem 2.6 of \cite{KaiErg}, any quasi-product $G$-invariant Radon measure on the double boundary of a CAT$(-1)$ space which gives full measure to pairs of conical limit points of $G$ is ergodic. Thus, $\nu \times \nu$ is ergodic.
\end{proof}

\section{RAAGs, RACGs and graph products} \label{RAAGs}

Let $\Lambda$ be a finite simplicial (undirected) graph. Recall that the corresponding right-angled Artin group (RAAG) $A(\Lambda)$ is the group given by the presentation
\[
A(\Lambda) := \langle v \in V(\Lambda) : [v,w] =1 \iff (v,w) \in E(\Lambda) \rangle.
\]
The corresponding right-angled Coxeter group (RACG) $C(\Lambda)$ is the group obtained from $A(\Lambda)$ by adding the relators $v^2 =1$ for each $v \in V(\Lambda)$.
In each case, $S = \{v^{\pm1}: v \in V(\Lambda)\}$ is called the set of \emph{standard} (or \emph{vertex}) generators of the group. 

In greater generality, let $\Lambda$ be a finite simplicial graph, and for each vertex $v$ of $\Lambda$ let us pick a finitely generated group $G_v$, which we call \emph{vertex group}.
Then we define the \emph{graph product} 
$$G(\Lambda) := \langle g \in G_v \ : \ [g,h] = 1 \iff g \in G_v, h \in G_w \textup{ and }(v, w) \in E(\Lambda) \rangle$$
as the group generated by the vertex groups $G_v$ with the relation that two vertex groups commute if and only if the corresponding vertices 
are joined by an edge. Clearly, RAAGs are special cases of graph products when $G_x = \mathbb{Z}$ for all $x$, and RACGs are graph products with $G_x = \mathbb{Z}/2\mathbb{Z}$. Graph products were first introduced by Green \cite{green1990graph} and have received much attention, see for example \cite{baik1993identity,chiswell1994growth,green1990graph,Hag,HerM,hsu1999linear,meier1996graph,meinert1995bieri,radcliffe2003rigidity}.

In this section, we are going to apply our counting techniques to graph products. 

\subsection{Geodesic combing for graph products} 
Let us call a group \emph{admissible} if it has a geodesic combing with respect to some finite generating set (i.e., in the language of the previous sections, if it 
has an admissible generating set). Recall that a recurrent component is nontrivial if it contains at least one closed path.
A component is terminal if there is no path exiting it. 
A graph structure is \emph{recurrent} if every vertex admits a directed path to every vertex other than the initial one.

Recall that, given a graph $\Lambda$, the \emph{opposite graph} is the graph $\Lambda^{op}$ with the same vertex set as $\Lambda$ and 
such that $(v, w) \in E(\Lambda^{op})$ if and only if $(v, w) \notin E(\Lambda)$.
We will assume that $\Lambda$ is \emph{anticonnected}, i.e. that the opposite graph $\Lambda^{op}$ is connected. 
This implies that $G(\Lambda)$ is not a direct product of graph products associated to subgraphs of $\Lambda$. 

\begin{proposition} \label{P:graph-combing}
Let $\Lambda$ be anticonnected, and choose for each vertex $x$ a group $G_x$ with a 
geodesic combining $(G_x, \Gamma_x)$
for the generating set $S_x$. 
Then the graph product $G(\Lambda)$ with the generating set $S = \cup_x S_x$ admits a 
geodesic combing which is recurrent. 
\end{proposition}

We call the generating set $S$ in Proposition \ref{P:graph-combing}, the \emph{standard} generating set for $G(\Lambda)$. Note that this agrees with the standard vertex generators for the special case of right-angled Artin and Coxeter groups.

The proof of Proposition \ref{P:graph-combing} will provide an explicit construction of a recurrent graph structure for $G(\Lambda)$ with the standard generators. Hermiller-Meier \cite{HerM} provided a construction of a geodesic combing, which is however not recurrent. 
In the next few lemmas, we will show that if $\Lambda$ is anticonnected we can 
modify their construction in order to make it recurrent. Of course it is necessary to assume that $\Lambda$ is anticonnected, as the counting theorems fail 
for RAAGs which decompose as direct products (see Example \ref{ex:nongen}).

Let us first review the construction of \cite{HerM}.
First introduce a total ordering on the vertices of $\Lambda$ such that the first two vertices in the ordering are not adjacent in $\Lambda$.
Each vertex of $\Lambda$ will be labeled by a capital letter $A, B, \dots$.

Then for each pair of vertices $(I, J)$ such that $I$ and $J$ are not adjacent in $\Lambda$ and with $I > J$ one constructs the $(I, J)$\emph{-admissible tree} in the following way. 
An $(I,J)$\emph{-admissible word} is a finite sequence $IJ K_1 K_2 \dots K_r$ (with $r \geq 0)$ such that:
\begin{enumerate}
\item 
$$J < K_1 < K_2 < \dots < K_r$$
and 
\item
if $K_i \leq I$ for some $i \leq r$, then $K_i$ is not adjacent to at least one vertex among $I, J, K_1, \dots, K_{i-1}$.
\end{enumerate}
Given $(I, J)$, the $(I, J)$-admissible tree is the finite directed tree, whose vertices are labeled by letters and whose paths spell exactly the $(I, J)$-admissible words.
In particular, such a tree will have $I$ as a root and there is only one edge coming out of this vertex, with endpoint $J$. Here, and in what follows, a directed edge always has the same label as its terminal vertex.

Moreover, \cite{HerM} define the \emph{header graph} (the terminology is ours) as the graph with one vertex for each letter, and an edge $A \to B$ if and only if $A< B$. 

Finally, they construct the graph structure for $C(\Lambda)$ (the corresponding RACG) as follows.
Consider the union of an initial vertex $v_0$, the header graph and all $(I, J)$-admissible trees. First, one identifies the vertex $I$ of the header graph with 
the root of the $(I, J)$-admissible tree for each possible $J$. Then, one adds one edge from $v_0$ to each vertex of the header graph, and 
if $A > B$ and $A, B$ are not adjacent, one joins by a directed edge each vertex labeled $A$ in the union of the $(I, J)$-admissible trees with the $B$ vertex in the 
$(A, B)$-admissible tree.
As shown by Hermiller-Meier, this graph $\G$ gives a bijective, geodesic graph structure for $C(\Lambda)$ \cite[Section 5 and Proposition 3.3]{HerM}. In fact, they show that $\G$ recognizes the geodesic language of normal forms with respect to the ordering of the vertices of $\Lambda$, but we will not need this stronger fact.

Let $C$ be the subgraph of $\G$ obtained by removing all vertices in the header graph and the initial vertex. 
That is, $C$ is the subgraph induced on the vertices on all $(I,J)$--admissible trees, excluding the initial vertex of the tree (which is labeled $I$).

\begin{lemma} \label{L:irred}
If $\Lambda$ is anticonnected, then the graph $C$ is irreducible, i.e. there is one directed path from each vertex to any other vertex.
Hence, $\G$ has a unique nontrivial recurrent component $C$ and this component is terminal. 
\end{lemma}

\begin{proof}
 We will show that $C$ is indeed irreducible. This will suffice since the header graph has no directed loops (it only has directed edges which increase in the ordering) and there are no edges leaving $C$ by construction.

Since the (unique) type $J$ vertex in the $(I,J)$--admissible tree has a directed path to each of its vertices and each vertex is in some $(I,J)$--admissible tree, it suffices to show that from any vertex of $C$ we can reach the type $J$ vertex of any $(I,J)$--admissible tree. Hence, fix some vertex $v$ of $\G$ and $I,J$, which are vertices of $\Lambda$.

Here is a main point: Any type $I$ vertex of any admissible tree is joined to the type $J$ vertex of the $(I,J)$--admissible tree. Hence, it suffices to get from $v$ to any type $I$ vertex of any admissible tree. To do this let $X$ be the type of $v$.

Fix a path $X = X_0, X_1, \ldots X_n = I$ in the complement graph of $\Lambda$. (Here we use that $\Lambda$ is anticonnected.) That is, $X_i$ and $X_{i+1}$ are not adjacent in $\Lambda$. We have to get from $v$ to any type $I$ vertex. We do this inductively as follows: We have either $X<X_1$ or $X>X_1$. In the first case, there is a type $X_1$ vertex $v_1$ in the admissible tree containing $v$ along with a directed path $v \to v_1$. (Since no consecutive pair in our fixed path are adjacent in $\Lambda$, condition (2) above holds automatically.) In the second case, there is an edge from $v$ to the unique $X_1$ vertex of the $(X,X_1)$--admissible tree, call this vertex $v_1$. In either case we get a directed path $v \to v_1$, were $v_1$ is a type $X_1$ vertex. 

We now repeat this argument to produce a path from $v_1 \to v_2$, where $v_2$ is a type $X_2$ vertex. Continuing in this manner, we produce $v \to v_1 \to \ldots \to v_n$, were $v_n$ has type $I$. Since $v_n$ then has a directed edge to the type $J$ vertex of the $(I,J)$--admissible tree (as discussed above), this completes the proof,
\end{proof}

We now know that the union of the admissible trees (excluding the initial vertices) is an irreducible graph. However, the header graph by construction is not irreducible. However, 
in the following lemma we observe that all words we can spell in the header graph can also be spelled in one of the admissible trees. Hence, we can modify $\G$ 
(essentially, by removing the header graph) in order to get a recurrent graph $\G^r$ which recognizes the same language as $\G$.

\begin{lemma}\label{lem:same_lang}
If $\Lambda$ is anticonnected, there exists a recurrent graph $\G^r$ which recognizes the same language as $\G$. 
\end{lemma}

\begin{proof}
Assume that $\Lambda$ is anticonnected and that its vertices are ordered so that the first two vertices $A, B$ do not commute (i.e. they are not adjacent). 
We modify $\G$ so that the resulting graph $\G^r$ still recognizes the same language as $\G$, and it is recurrent.

The modification is simple and requires only one observation: we note that any strictly increasing sequence $X_1 \ldots X_r$ can be spelled in the $(B,A)$--admissible tree, starting from some vertex. In fact, if $X_1 = A$ then $B X_1 \ldots X_r$ is $(B, A)$-admissible, since the only required condition is that whenever $X_i \leq B$ the vertex $X_i$ is not adjacent to some $X_l$ with $l < i$. However, the only two letters not greater than $B$ are $A, B$, and $A$ and $B$ are not adjacent by construction. 
Similarly,  if $X_1 \neq A$ then $BA X_1 \ldots X_r$ is $(B, A)$-admissible.

Thus, the new graph $\G^r$ is given by removing the header graph and joining the initial vertex $v_0$ to each vertex of the $(B, A)$-admissible tree. Any word which is recognized by $\G$ is made of an increasing word followed by a word spelled in the union of the admissible trees. In $\G^r$, such a word is spelled by spelling the increasing sequence in the $(B, A)$-admissible tree, 
and the second part as before. This proves the claim.
\end{proof}

The graph $\G^r$ is a recurrent graph which by Lemma \ref{lem:same_lang} gives a bijective, geodesic graph structure the right-angled Coxeter group $C(\Lambda)$. We now modify the construction 
to produce a geodesic combing for each graph product $G(\Lambda)$.

Let $\Gamma_I$ be the graph structure of the vertex group $G_I$, let $v_{0, I}$ be the initial vertex of $G_I$,
let $s_{1, I}, \dots, s_{k, I}$ be the labels of the edges going out of $v_{0, I}$, and let $v_{1, I}, \dots, v_{I, k}$ be the targets of these edges, respectively.
Moreover, let $\Gamma'_I$ be the subgraph of $\Gamma_I$ given by removing the initial vertex.

To construct the graph structure for $G(\Lambda)$, let us consider the disjoint union of a vertex $\widetilde{v}_0$, which will serve as initial vertex, 
and a copy of $\Gamma_I'$ for each vertex $v$ of type $I$ in $\G^r$.
Moreover, for any edge in $\G^r$ of type $I \to J$ let us connect each vertex of the corresponding $\Gamma'_I$ with the vertices $v_{1, J}, \dots, v_{k, J}$ of 
the corresponding $\Gamma'_J$ with edges labeled, respectively, $s_{1, J}, \dots, s_{k, J}$. 
Finally, for each edge from $v_0$ in $\G^r$ to some other vertex of type $I$, let us connect the new initial vertex $\widetilde{v}_0$ with vertices $v_{1, J}, \dots, v_{k, J}$ of $\Gamma'_J$ with edges labeled, respectively, $s_{1, J}, \dots, s_{k, J}$. 

This new graph $\Gamma_G$ gives a bijective, geodesic structure for $G(\Lambda)$ with respect to the standard generators. This follows since, by construction, $\Gamma_G$ parameterizes the same language of geodesic normal forms for $G(\Lambda)$ given in \cite{HerM}. Moreover, since $\Gamma_G$ is modeled on the recurrent graph $\G^r$ one easily sees that $\Gamma_G$ is itself recurrent.
This completes the proof of Proposition \ref{P:graph-combing}.
 
\begin{corollary}
Let $G(\Lambda)$ be a graph product of admissible groups, which does not decompose as a direct product.
Then there exists a thick graph structure for its standard generating set. 
\end{corollary}

\begin{proof} 
From Proposition \ref{P:rec-thick}, the graph structure given by the above proposition is thick since $\Gamma_G$ is recurrent.
\end{proof}

As a consequence of thickness, we are ready to establish the following counting result for loxodromics. 

\begin{theorem} \label{th:RAAGs}
Let $G$ be an infinite graph product of admissible groups 
which does not decompose as a product of infinite groups, and let $S$ be its set of standard (vertex) generators. 
Then for any nonelementary action $G \curvearrowright X$ on a separable hyperbolic space $X$, the set of loxodromics for the action is generic with respect to $S$, i.e.
\[
\frac{ \#\{g \in G : |g|_S \le n \text{ and } g \text{ is } X - \text{loxodromic}\}}{\#\{g\in G : |g|_S\le n \}} \longrightarrow 1,
\]
as $n \to \infty$.
\end{theorem}

\section{Exact exponential growth for RAAGs and RACGs}

We conclude by proving a fine estimate on the number of elements in a ball for RAAGs and RACGs.

\begin{theorem} \label{T:exact-exp}
Let $G$ be a right-angled Artin or Coxeter group which is not virtually cyclic and does not decompose as a product of infinite groups, and let us consider 
$S$ its standard generating set. Then there exists constants $\lambda > 1$, $c > 0$ such that the following limit exists: 
\begin{equation} \label{E:exact-exp}
\lim_{n \to \infty} \frac{\#\{ g \in G \ : \ |g|_S = n\}}{\lambda^n} = c.
\end{equation}
\end{theorem}

We say that a group which satisfies \eqref{E:exact-exp} has \emph{exact exponential growth}. Let us remark that such a property 
is \emph{not} invariant with respect to quasi-isometries of the metric, and hence it depends very carefully on the generating set.

In fact, Theorem \ref{T:exact-exp} will follow immediately from the following theorem for general graph products. 

\begin{theorem} \label{T:exact-graph}
Let $G(\Lambda)$ be a graph product of admissible groups, and assume that $\Lambda$ is anticonnected (so that the group does not split trivially as a product) and 
has at least $3$ vertices. Then $G(\Lambda)$ has exact exponential growth. 
\end{theorem}

Note that it makes sense to assume that the number of vertices is at least $3$. In fact, if $n = 1$ then $G(\Lambda)$ can be any group with a geodesic combing, 
while if $n = 2$ then $G(\Lambda)$ can be the free product of any two admissible groups. In particular, if it is a RAAG then it must be the free group on $2$ generators, 
which has exact exponential growth, and if it is a RACG it must be $\mathbb{Z}/2\mathbb{Z} \star \mathbb{Z}/2\mathbb{Z}$, which is virtually cyclic.

Let us remark that the growth function for graph products has been worked out by Chiswell \cite{chiswell1994growth} (see also \cite{Athreya-Prasad}); however, it does not seem obvious how to prove exact exponential growth by this method. 

Let us consider the recurrent graph $\G^r$ defined in the previous section, and denote as $\G^r_0 = \G^r \setminus \{ v_0 \}$. By the previous section, we know that $\G^r_0$ is irreducible. The final step in the proof of Theorem \ref{T:exact-graph} is the following lemma.

\begin{lemma} \label{L:aper}
If $\Lambda$ is anticonnected and has at least $3$ vertices, then the graph $\G^r_0$ is aperiodic. 
\end{lemma}

\begin{proof}
Let us assume, consistently with the previous section, that the vertices of $\Lambda$ are ordered. Let us call $A, B, C$ 
the three smallest vertices, with $A < B < C$, and assume that $A$, $B$ are not adjacent. 
Then let us observe that the sequences $BAC$ and $BABC$ are $(B, A)$ admissible, hence in the $(B, A)$-admissible tree there is a $Y$-shaped subtree 
with five vertices: one labeled $A$, two labeled $B$ (let us denote them $B_1$, $B_2$) and two labeled $C$ (let us denote them $C_1$, $C_2$) so that 
the paths in this subtree are $B_1 \to A \to C_1$ and $B_1 \to A \to B_2 \to C_2$. 
Now, since the graph is irreducible, there exists a path from $C_1$ to $A$; let us denote its vertices as 
$C_1 \to v_1 \to v_2 \dots \to v_k \to A$. Then by definition, the type of $v_1$ is smaller than $C$, and is not adjacent to $C$. Thus, by construction, there is also an edge 
from $C_2$ to $v_1$; hence, in the graph 
there are two loops: one loop is given by $A \to C_1 \to v_1\to  \dots \to v_k \to A$ and the other is $A \to B \to C_2 \to v_1 \to  \dots \to v_k \to A$. Since the lengths of these two closed paths differ by one, the greatest common divisor 
of the lengths of all paths is $1$, hence $\G^r_0$ is aperiodic. 
\end{proof}

Note that the statement is false if the number of vertices is $2$: indeed, then there is only one loop of length $2$, hence the period is $2$.

Now, let us consider a general graph product $G(\Lambda)$. By the previous section, by replacing vertices of $\G^r$ with graphs which recognize the geodesic 
combings of vertex group, we get a new graph $\Gamma_G$ which gives a geodesic combing for $G(\Lambda)$. By the previous Lemma we get: 

\begin{corollary}
If $\Lambda$ is anticonnected and has at least three vertices, then the graph $\Gamma_G' = \Gamma_G \setminus \{ v_0 \}$ is irreducible and aperiodic. 
\end{corollary}

\begin{proof}[Proof of Theorem \ref{T:exact-graph}]
Since the graph $\Gamma'_G$ is irreducible and aperiodic, then by the Perron-Frobenius theorem its adjacency matrix $A$ has a unique eigenvalue $\lambda > 1$ of maximum modulus,  
and that eigenvalue is real, positive, and simple. Moreover, the coordinates of the corresponding eigenvector are all positive. Finally, the sequence 
$\frac{A^n}{\lambda^n}$ converges to the projection to the eigenspace. In particular, none of the basis vectors is orthogonal to the eigenvector, hence for any $i,j$ there exists $c_{ij} > 0$ such that 
$$\lim_{n \to \infty} \frac{(A^n)_{ij}}{\lambda^n} = c_{ij}.$$ 
Now, each path of length $n$ from the initial vertex starts with an edge to the irreducible graph, hence 
$$\frac{\#S_n}{\lambda^n} = \sum_{v_0 \to v_i} \frac{\# S_{n-1}(v_i)}{\lambda^n} \to \sum_{v_0 \to v_i} \sum_j \frac{(A^{n-1})_{ij}}{\lambda^n} \to \sum_{v_0 \to v_i} \sum_j \frac{c_{ij}}{\lambda} = c > 0$$
which establishes exact exponential growth.
\end{proof}

\bibliography{counting2}

\newcommand{\etalchar}[1]{$^{#1}$}
\providecommand{\bysame}{\leavevmode\hbox to3em{\hrulefill}\thinspace}
\providecommand{\MR}{\relax\ifhmode\unskip\space\fi MR }
\providecommand{\MRhref}[2]{%
  \href{http://www.ams.org/mathscinet-getitem?mr=#1}{#2}
}
\providecommand{\href}[2]{#2}
\begin{thebibliography}{EPC{\etalchar{+}}92}

\bibitem[AC16]{AntolinCiobanu}
Yago Antolin and Laura Ciobanu, \emph{Finite generating sets of relatively
  hyperbolic groups and applications to geodesic languages}, Trans. Amer. Math.
  Soc. \textbf{368} (2016), no.~11, 7965--8010.

\bibitem[ACT15]{Arzhantseva-Cashen-Tao}
G.N. Arzhantseva, C.~Cashen, and J.~Tao, \emph{Growth tight actions}, Pacific
  J. Math. \textbf{278} (2015), no.~1, 1--49.

\bibitem[ADT]{ADT}
Tarik Aougab, Matthew~Gentry Durham, and Samuel~J. Taylor, \emph{Pulling back
  stability with applications to {O}ut({F}) and relatively hyperbolic groups},
  J. Lond. Math. Soc.

\bibitem[AL02]{arzhantseva2002growth}
Goulnara~N Arzhantseva and IG~Lysenok, \emph{Growth tightness for word
  hyperbolic groups}, Mathematische Zeitschrift \textbf{241} (2002), no.~3,
  597--611.

\bibitem[AO96]{ArzhantsevaOlshanskii}
G.~Arzhantseva and A.~Olshanskii, \emph{Generality of the class of groups in
  which subgroups with a lesser number of generators are free}, Mat. Zametki
  \textbf{59} (1996), 489--496.

\bibitem[AP14]{Athreya-Prasad}
J.~Athreya and A.~Prasad, \emph{Growth in right-angled groups and monoids},
  available at arXiv:1409.4142, 2014.

\bibitem[Arz98]{arzhantseva1998generic}
GN~Arzhantseva, \emph{Generic properties of finitely presented groups and
  {H}owson's theorem}, Communications in Algebra \textbf{26} (1998), no.~11,
  3783--3792.

\bibitem[BF14]{BF14}
Mladen Bestvina and Mark Feighn, \emph{Hyperbolicity of the complex of free
  factors}, Adv. Math. \textbf{256} (2014), 104--155. \MR{3177291}

\bibitem[BH09]{BH}
Martin~R. Bridson and Aandre Haefliger, \emph{Metric spaces of non-positive
  curvature}, vol. 319, Springer, 2009.

\bibitem[BHP93]{baik1993identity}
Young-Gheel Baik, Jim Howie, and SJ~Pride, \emph{The identity problem for graph
  products of groups}, Journal of Algebra \textbf{162} (1993), no.~1, 168--177.

\bibitem[BMR03]{borovik2003multiplicative}
Alexandre~V Borovik, Alexei~G Myasnikov, and Vladimir~N Remeslennikov,
  \emph{Multiplicative measures on free groups}, Internat. J. Algebra Comput.
  \textbf{13} (2003), no.~06, 705--731.

\bibitem[Bow08]{Bo}
Brian~H. Bowditch, \emph{Tight geodesics in the curve complex}, Invent. Math.
  \textbf{171} (2008), no.~2, 281--300.

\bibitem[Bow12]{bowditch2012relatively}
Brian~H Bowditch, \emph{Relatively hyperbolic groups}, International Journal of
  Algebra and Computation \textbf{22} (2012), no.~03, 1250016.

\bibitem[Cal13]{calegari2013ergodic}
Danny Calegari, \emph{The ergodic theory of hyperbolic groups}, Geometry and
  topology down under, Contemp. Math \textbf{597} (2013), 15--52.

\bibitem[Can84]{cannon1984combinatorial}
James~W Cannon, \emph{The combinatorial structure of cocompact discrete
  hyperbolic groups}, Geometriae Dedicata \textbf{16} (1984), no.~2, 123--148.

\bibitem[CF10]{calegari2010combable}
Danny Calegari and Koji Fujiwara, \emph{Combable functions, quasimorphisms, and
  the central limit theorem}, Ergodic Theory and Dynamical Systems \textbf{30}
  (2010), no.~05, 1343--1369.

\bibitem[Cha95]{champetier1995proprietes}
Christophe Champetier, \emph{Propri{\'e}t{\'e}s statistiques des groupes de
  pr{\'e}sentation finie}, Adv. Math. \textbf{116} (1995), no.~2, 197--262.

\bibitem[Chi94]{chiswell1994growth}
Ian~M Chiswell, \emph{The growth series of a graph product}, Bulletin of the
  London Mathematical Society \textbf{26} (1994), no.~3, 268--272.

\bibitem[CM07]{connell2007harmonicity}
Chris Connell and Roman Muchnik, \emph{Harmonicity of quasiconformal measures
  and {P}oisson boundaries of hyperbolic spaces}, Geom. Funct. Anal.
  \textbf{17} (2007), no.~3, 707--769.

\bibitem[CM15]{calegari2015statistics}
Danny Calegari and Joseph Maher, \emph{Statistics and compression of scl},
  Ergodic Theory and Dynamical Systems \textbf{35} (2015), no.~01, 64--110.

\bibitem[DMS10]{DMS}
Cornelia Dru{\c{t}}u, Shahar Mozes, and Mark Sapir, \emph{Divergence in
  lattices in semisimple {L}ie groups and graphs of groups}, Transactions of
  the American Mathematical Society \textbf{362} (2010), no.~5, 2451--2505.

\bibitem[DSU14]{das2014geometry}
Tushar Das, David Simmons, and Mariusz Urba{\'n}ski, \emph{Geometry and
  dynamics in {G}romov hyperbolic metric spaces: With an emphasis on non-proper
  settings}, arXiv preprint arXiv:1409.2155 (2014).

\bibitem[DT06]{DunfieldThurston}
Nathan Dunfield and William~P Thurston, \emph{Finite covers of random
  3-manifolds}, Invent. Math. \textbf{166} (2006), no.~3, 457--521.

\bibitem[EPC{\etalchar{+}}92]{epstein1992word}
David Epstein, Mike~S Paterson, James~W Cannon, Derek~F Holt, Silvio~V Levy,
  and William~P Thurston, \emph{Word processing in groups}, AK Peters, Ltd.,
  1992.

\bibitem[Far98]{farb1998relatively}
Benson Farb, \emph{Relatively hyperbolic groups}, Geom. Funct. Anal. \textbf{8}
  (1998), no.~5, 810--840.

\bibitem[Fur02]{Furman}
A.~Furman, \emph{Coarse-geometric perspective on negatively curved manifolds
  and groups}, Rigidity in Dynamics and Geometry (M.~Burger and A.~Iozzi,
  eds.), Springer, 2002, pp.~149--166.

\bibitem[GdlH90]{GhysdelaHarpe}
{\'E}.~Ghys and P.~de~la Harpe (eds.), \emph{Sur les groupes hyperboliques
  d'apr\`es {M}ikhael {G}romov}, Progress in Mathematics, vol.~83, Birkh\"auser
  Boston, Inc., Boston, MA, 1990, Papers from the Swiss Seminar on Hyperbolic
  Groups held in Bern, 1988. \MR{1086648 (92f:53050)}

\bibitem[GM08]{groves2008dehn}
Daniel Groves and Jason~Fox Manning, \emph{Dehn filling in relatively
  hyperbolic groups}, Israel Journal of Mathematics \textbf{168} (2008), no.~1,
  317--429.

\bibitem[GMM15]{gouezel2015entropy}
S{\'e}bastien Gou{\"e}zel, Fr{\'e}d{\'e}ric Math{\'e}us, and Fran{\c{c}}ois
  Maucourant, \emph{Entropy and drift in word hyperbolic groups}, arXiv
  preprint arXiv:1501.05082 (2015).

\bibitem[Gre90]{green1990graph}
Elisabeth~Ruth Green, \emph{Graph products of groups}, Ph.D. thesis, University
  of Leeds, 1990.

\bibitem[Gro87]{Gromov}
Mikhael Gromov, \emph{Hyperbolic groups}, Springer, 1987.

\bibitem[Gro93]{gromov1993geometric}
M~Gromov, \emph{Asymptotic invariants of infinite groups}, Geometric group
  theory, vol. 2, vol. 182, Cambridge University Press, 1993.

\bibitem[Gro03]{gromov2003random}
Mikhail Gromov, \emph{Random walk in random groups}, Geometric and Functional
  Analysis \textbf{13} (2003), no.~1, 73--146.

\bibitem[GTT16]{GTT}
I.~Gekhtman, S.~Taylor, and G.~Tiozzo, \emph{Counting loxodromics for
  hyperbolic actions}, available at arXiv:1605.02103, 2016.

\bibitem[Hag08]{Hag}
F.~Haglund, \emph{Finite index subgroups of graph products}, Geometriae
  Dedicata \textbf{135} (2008), no.~1, 167--209.

\bibitem[Hag14]{hagen2014weak}
Mark~F Hagen, \emph{Weak hyperbolicity of cube complexes and quasi-arboreal
  groups}, Journal of Topology \textbf{7} (2014), no.~2, 385--418.

\bibitem[HM95]{HerM}
S.~Hermiller and J.~Meier, \emph{Algorithms and geometry for graph products of
  groups}, Journal of Algebra \textbf{171} (1995), no.~1, 230--257.

\bibitem[HM13]{handel2013free}
Michael Handel and Lee Mosher, \emph{The free splitting complex of a free
  group, i hyperbolicity}, Geometry \& Topology \textbf{17} (2013), no.~3,
  1581--1670.

\bibitem[HRR17]{GroupsLangAuto}
Derek~F. Holt, Sarah Rees, and Claas~E. R\"over, \emph{Groups, languages and
  automata}, London Mathematical Society Student Texts, vol.~88, Cambridge
  University Press, Cambridge, 2017.

\bibitem[Hru10]{hruska}
G.~Hruska, \emph{Relative hyperbolicity and relative quasiconvexity for
  countable groups}, Algebr. Geom. Topol \textbf{10} (2010), 1807--1856.

\bibitem[HW{\etalchar{+}}99]{hsu1999linear}
Tim Hsu, Daniel~T Wise, et~al., \emph{On linear and residual properties of
  graph products.}, The Michigan Mathematical Journal \textbf{46} (1999),
  no.~2, 251--259.

\bibitem[Kai94]{KaiErg}
Vadim Kaimanovich, \emph{Ergodicity of harmonic invariant measures for the
  geodesic flow on hyperbolic spaces}, J. Reine Angew. Math. \textbf{455}
  (1994), 57--104.

\bibitem[KB02]{kap_boundaries}
Ilya Kapovich and Nadia Benakli, \emph{Boundaries of hyperbolic groups},
  Combinatorial and geometric group theory (New York, 2000/Hoboken, NJ, 2001)
  \textbf{296} (2002), 39--93.

\bibitem[KK14]{KK3}
Sang-Hyun Kim and Thomas Koberda, \emph{The geometry of the curve graph of a
  right-angled {A}rtin group}, International Journal of Algebra and Computation
  \textbf{24} (2014), no.~02, 121--169.

\bibitem[KMSS03]{kapovich2003generic}
Ilya Kapovich, Alexei Myasnikov, Paul Schupp, and Vladimir Shpilrain,
  \emph{Generic-case complexity, decision problems in group theory, and random
  walks}, J. Algebra \textbf{264} (2003), no.~2, 665--694.

\bibitem[KRSS07]{kapovich2007densities}
Ilya Kapovich, Igor Rivin, Paul Schupp, and Vladimir Shpilrain, \emph{Densities
  in free groups and $\mathbb{Z}^n$, visible points and test elements}, Math.
  Res. Lett \textbf{14} (2007), no.~2.

\bibitem[KS05]{kapovich2005genericity}
Ilya Kapovich and Paul Schupp, \emph{Genericity, the
  {A}rzhantseva-{O}l'shanskii method and the isomorphism problem for
  one-relator groups}, Math. Ann. \textbf{331} (2005), no.~1, 1--19.

\bibitem[Mah11]{Maher}
Joseph Maher, \emph{Random walks on the mapping class group}, Duke Math. J.
  \textbf{156} (2011), no.~3, 429--468.

\bibitem[Mei95]{meinert1995bieri}
Holger Meinert, \emph{The bieri-neumann-strebel invariant for graph products of
  groups}, Journal of Pure and Applied Algebra \textbf{103} (1995), no.~2,
  205--210.

\bibitem[Mei96]{meier1996graph}
John Meier, \emph{When is the graph product of hyperbolic groups hyperbolic?},
  Geometriae Dedicata \textbf{61} (1996), no.~1, 29--41.

\bibitem[MM99]{MM1}
Howard~A. Masur and Yair~N. Minsky, \emph{Geometry of the complex of curves.
  {I}. {H}yperbolicity}, Invent. Math. \textbf{138} (1999), no.~1, 103--149.

\bibitem[MS14]{mathieu2014deviation}
P~Mathieu and A~Sisto, \emph{Deviation inequalities for random walks}, arXiv
  preprint arXiv:1411.7865 (2014).

\bibitem[MT14]{MaherTiozzo}
Joseph Maher and Giulio Tiozzo, \emph{Random walks on weakly hyperbolic
  groups}, To appear in J. Reine Angew. Math. (2014).

\bibitem[Neu54]{Neumann}
Bernhard~H. Neumann, \emph{Groups covered by finitely many cosets}, Publ. Math.
  Debrecen \textbf{3} (1954).

\bibitem[NS95]{Neumann-Shapiro-kleinian}
Walter Neumann and Michael Shapiro, \emph{Automatic structures, rational growth
  and geometrically finite hyperbolic groups}, Invent. Math. \textbf{120}
  (1995), 259--287.

\bibitem[Ol'92]{ol1992almost}
A~Yu Ol'shanskii, \emph{Almost every group is hyperbolic}, Internat. J. Algebra
  Comput. \textbf{2} (1992), no.~01, 1--17.

\bibitem[Osi06]{osin2006relatively}
Denis~V Osin, \emph{Relatively hyperbolic groups: intrinsic geometry, algebraic
  properties, and algorithmic problems}, vol. 843, American Mathematical Soc.,
  2006.

\bibitem[Osi15]{osin2015acylindrically}
Denis Osin, \emph{Acylindrically hyperbolic groups}, Trans. Amer. Math. Soc.
  (2015).

\bibitem[Rad03]{radcliffe2003rigidity}
David~G Radcliffe, \emph{Rigidity of graph products of groups}, Algebraic \&
  Geometric Topology \textbf{3} (2003), no.~2, 1079--1088.

\bibitem[Riv08]{rivin2008walks}
Igor Rivin, \emph{Walks on groups, counting reducible matrices, polynomials,
  and surface and free group automorphisms}, Duke Math. J. \textbf{142} (2008),
  no.~2, 353--379.

\bibitem[Sis11]{sisto2011contracting}
Alessandro Sisto, \emph{Contracting elements and random walks}, J. Reine Angew.
  Math. (2011).

\bibitem[TT15]{TT}
Samuel~J. Taylor and Giulio Tiozzo, \emph{Random extensions of free groups and
  surface groups are hyperbolic}, Int. Math. Res. Not. (2015).

\bibitem[Wie14]{wiest2014genericity}
Bert Wiest, \emph{On the genericity of loxodromic actions}, arXiv preprint
  arXiv:1406.7041 (2014).

\bibitem[Yan13]{YangPS}
Wenyuan Yang, \emph{{P}atterson-{S}ullivan measures and growth of relatively
  hyperbolic groups}, available at arXiv:1308.6326, 2013.

\bibitem[Yan16]{Wenyuan-scc}
\bysame, \emph{Statistically convex-cocompact actions of groups with
  contracting elements}, available at arXiv:1612.03648, 2016.

\end{thebibliography}
\bibliographystyle{amsalpha}

\end{document}